\documentclass[11pt]{article}
\usepackage[utf8]{inputenc}

\usepackage[numbers,sort&compress]{natbib}
\usepackage{authblk}

\usepackage{algorithm}
\usepackage{algorithmic}

\usepackage{pgf}%figure
\usepackage{tikz}
\usetikzlibrary{arrows, decorations.pathmorphing, backgrounds, positioning, fit, petri, automata}
\definecolor{yellow1}{rgb}{1,0.8,0.2}

\usepackage{makecell}

\usepackage{amsmath,amsthm,amsfonts,amssymb,amsbsy,amsopn,amstext}
\usepackage{graphicx,color}
\usepackage{mathbbol,mathrsfs}
\usepackage{float,ccaption}
\usepackage{indentfirst}
\usepackage{appendix}
\usepackage{subfigure}

\usepackage{emptypage}

\usepackage{booktabs}
\usepackage{threeparttable}
\usepackage{multirow}
\usepackage{diagbox}
\renewcommand{\arraystretch}{1.5} %控制行高

 \newtheorem{thm}{Theorem}
 
 \newtheorem{lem}{Lemma}
 \newtheorem{prop}{Proposition}

 \newtheorem{ass}{Assumption}
\newcommand{\e}{\ensuremath{\mathrm{e}}} %Euler????,???pi??\piup
\renewcommand{\P}{\ensuremath{\mathrm{P}}} %????P
\newcommand{\define}{:=}%\newcommand{\define}{\stackrel{\triangle}{=}}
\newcommand{\normm}[1]{{\left\vert\kern-0.25ex\left\vert\kern-0.25ex\left\vert #1
		\right\vert\kern-0.25ex\right\vert\kern-0.25ex\right\vert}}

\DeclareMathOperator*{\diag}{diag}

\usepackage{geometry}
\begin{document}
\title{Efficient  Gradient Tracking Algorithms for Distributed Optimization Problems with Inexact Communication}

\author[a]{Shengchao Zhao}
\author[b]{Yongchao Liu \thanks{CONTACT Yongchao Liu.}}

\affil[a]{School of Mathematics, China University of Mining and Technology, Xuzhou, China\\zhaosc@cumt.edu.cn (Shengchao Zhao)}
\affil[b]{School of Mathematical Sciences, Dalian University of Technology, Dalian, China\\ lyc@dlut.edu.cn (Yongchao Liu)}

\date{}
\maketitle

\noindent{\bf Abstract.} 
Distributed optimization problems usually face   inexact communication issues induced by channel noise, communication quantization or differential privacy protection.
Most existing algorithms need a two-timescale setting of the stepsize of gradient descent  and  the parameter of noise suppression to ensure the convergence to the optimal solution. In this paper, we propose two single-timescale algorithms,  VRA-DGT and VRA-DSGT\footnote{In contexts involving both algorithms, they are jointly termed VRA-D(S)GT.}, for  distributed deterministic and stochastic optimization problems with inexact communication respectively.
VRA-DGT integrates the Variance-Reduced Aggregation (VRA) mechanism with the distributed gradient tracking framework,  which achieves the convergence rate of $\mathcal{O}\left(k^{-1}\right)$ in the mean square sense and  $\mathcal{O}\left(\frac{\ln(k+1)}{k^b}\right)$, $\forall b\in(0.5,1)$ in the almost sure sense when the objective function is strongly convex and smooth.
For stochastic optimization problems,  VRA–DSGT,  where a hybrid variance-reduced technique has been introduced in VRA-DGT, maintains  the convergence rate of $\mathcal{O}\left(k^{-1}\right)$ in the mean square sense and  $\mathcal{O}\left(\frac{\ln(k+1)}{k^b}\right)$, $\forall b\in(0.5,1)$ in the almost sure sense.
Simulated experiments on a logistic regression problem with real-world data verify the effectiveness of the proposed algorithms.	

\noindent\textbf{Key words.} distributed optimization, inexact communication, distributed gradient tracking

\section{Introduction}	

Distributed learning, whose implementation relies on communication over wireless channels, has garnered increasing interest owing to various applications that involve local data \cite{Saha2022Noi,Ayaz2019IoT,Samara2020vehic}. 
In many practical scenarios, the communication is performed under imperfect
environments,  for instance, wireless channels are subject to noise and fading \cite{Dasarathan2015robust,Kar2009imc,Saha2022Noi,Hong2021FL,Shah2023FL}; messages exchanged among agents are compressed to reduce communication costs \cite{Rabbat2004sensor,Koloskova2019spare};  random noise is added to transmitted information to protect privacy \cite{Wang2022Tailoring,Wang2017DPL}. In this paper, we consider the   distributed optimization problem
\begin{equation}\label{model}
	\min_{x\in\mathbb{R}^d} f(x)\define\frac{1}{n}\sum_{j=1}^n f_j(x),
\end{equation}
over the imperfect-communication network,  where $x$ is the global decision variable, $f_j(x)$ is the local cost
function specific to agent  $j\in\{1,2,\cdots,n\}$.  When the local cost function  $f_j(x):=\mathbb{E}\left[F_j(x;\xi_j)\right]$, where $\mathbb{E}\left[\cdot\right]$ represents the expectation taken over the probability space of random variable $\xi_j$, and $F_j(\cdot;\xi_j):\mathbb{R}^d\rightarrow \mathbb{R}$ is a measurable  function, problem (\ref{model}) denotes  the distributed stochastic optimization problem.  Given that each agent only has access to its local cost function, the goal of problem (1) is to find an optimal and consensual solution through the exchange of information over the imperfect-communication network.

A series of existing algorithmic schemes have been further developed to address distributed optimization problems with inexact communication, such as,  distributed dual averaging  \cite{YUAN2012DDA}, primal-dual  \cite{Lei2018}, Newton
tracking \cite{Pan2024Second}, gradient descent  \cite{Sri2011async} and gradient tracking \cite{pu2020robust}.
Among them,   Distributed Gradient Descent (DGD) and Distributed Gradient Tracking (DGT) are two popular algorithmic schemes. For the DGD scheme, Srivastava and Nedic \cite{Sri2011async} present a robust distributed stochastic (sub)gradient algorithm under general noisy communication environments, and establish the  almost sure convergence of the  algorithm for convex optimization problems. Reisizadeh et al. \cite{Reisizadeh2023Vary,Reisizadeh2023Dimix}
extend this robust algorithm to distributed nonconvex and strongly convex optimization problems, and achieve the convergence rate of  $\mathcal{O}(k^{-1/3+\delta})$ and $\mathcal{O}(k^{-1/2})$ respectively, where $k$ is the number of iteration and $\delta>0$ could be arbitrarily close to 0. In the context of   communication quantization,    Doan et al. \cite{Doan2021Quanti} develop a projection distributed gradient descent algorithm and establish the almost sure convergence of the proposed algorithm for constrained convex optimization problems. By introducing a novel analysis,  Vasconcelos et al. \cite{Vascon2021Improv} provide an improved convergence rate of $\mathcal{O}\left(\log^2 (k)/k^{1/2}\right)$  for the algorithm proposed in \cite{Doan2021Quanti}. On the other hand,  Reisizadeh et al. \cite{Reis2019Exact} consider the 
distributed strongly convex optimization problems with communication quantization and
propose a communication-efficient DGD algorithm by using a random quantizer, which achieves a convergence rate of $\mathcal{O}\left(k^{-1/2+\delta}\right)$.
More recently, Iakovidou and Wei \cite{Iakovi2023S-NEAR} propose a distributed stochastic gradient descent algorithm by utilizing the error correlation and nested consensus mechanisms to tolerate communication quantization, which converges to a neighborhood of the optimal solution with a linear rate.

While DGD based algorithms are simple to implement, they are sensitive to heterogeneous data/functions across agents \cite{koloskova2020Unified}. Noting that DGT algorithms can mitigate the impact of heterogeneous function through adding a vector to track the global gradient, Shah and Bollapragada \cite{shah2023SGT} propose a robust variant of the DGT algorithm, termed IC-GT, for distributed stochastic optimization problems, which suppress the information-sharing noise by incorporating small factors into the weight matrices of conventional DGT algorithm. When the objective function is strongly convex, IC-GT achieves a complexity of $\tilde{\mathcal{O}}(\epsilon^{-1})$ with stepsize $\alpha=\epsilon$.  However, IC-GT inherits from conventional DGT the issue of noise accumulation during gradient tracking 
To handle the noise accumulation issue, Pu \cite{pu2020robust} proposes the Robust Push-Pull  algorithm by tracking the cumulative global gradient rather than the global gradient. Under constant stepsize, Robust Push-Pull  algorithm converges linearly to a neighborhood of the optimal solution for strongly convex optimization problems.
To ensure  convergence to the optimal solution, Wang and Basar \cite{wangGT2022} propose a new robust DGT based algorithm by reconstructing
the cumulative gradient tracking to accommodate the decaying
factors, where the almost sure convergence to the optimal solution is obtained. Zhao et al. \cite{zhao2023VRA} present the VRA-GT algorithm by incorporating a Variance-Reduced Aggregation (VRA) mechanism into the cumulative global gradient tracking process of the R-Push-Pull algorithm, which achieves  the convergence rate of  $\mathcal{O}(k^{-1+\delta})$ for strongly convex optimization problems.

Indeed, the DGD and DGT based algorithms mentioned above  need a two-timescale condition on the stepsize of gradient descent and the parameter of noise suppression to ensure the convergence to the optimal solution. Specifically,  the parameter $\gamma_k$, which suppresses information-sharing noise, and the stepsize $\alpha_k$ of the gradient descent satisfy $\lim_{k\rightarrow\infty}\frac{\alpha_k }{\gamma_k }=0$. The condition implies that $\gamma_k$ varies significantly on a ``fast" timescale, while $\alpha_k $ changes slightly on a ``slow" timescale \cite{Doan2021Quanti}. Under the two-timescale condition $\lim_{k\rightarrow\infty}\frac{\alpha_k }{\gamma_k }=0$ and the condition $\sum_{k=1}^\infty \frac{\alpha_k^2}{\gamma_k} < \infty$ \cite{zhao2023VRA,wangGT2022,Doan2021Quanti}, the step size $\alpha_k$ must be sufficiently small, which results in slower convergence rates, that is,  $\mathcal{O}\left(k^{-1/2}\right)$ for DGD algorithms and $\mathcal{O}\left(k^{-1+\delta}\right)$ for DGT algorithms.
This motivates us to investigate single-timescale algorithms for distributed deterministic and stochastic optimization problems with inexact communication. We present the Variance-Reduced Aggregation-based Distributed Gradient Tracking (VRA-DGT) algorithm for distributed deterministic optimization problems and the Variance-Reduced Aggregation-based Distributed Stochastic Gradient Tracking (VRA-DSGT) algorithm for distributed stochastic optimization problems. As far as we are concerned, the contributions of the paper can be summarized as follows.
\begin{itemize}
	\item[(i)]For distributed  deterministic optimization problems with inexact communication, VRA-DGT integrates the VRA mechanism into both the decision variable updating and cumulative gradient tracking processes of R-Push-Pull \cite{pu2020robust}. In the decision variable updating process, VRA suppresses the noise generated during the sharing of decision variables among agents, which effectively alleviates the impact of communication error on decision variable updating. Similarly, in the cumulative gradient tracking process, VRA suppresses the noise that occurs when gradient trackers are shared between agents, which enables the dynamic consensus framework of gradient tracking to be preserved even in inexact communication environments. Consequently, VRA-DGT allows the parameter $\gamma_k$ of noise suppression and the stepsize $\alpha_k$ of gradient descent to satisfy the condition $\lim_{k\rightarrow\infty}\frac{\alpha_k }{\gamma_k }=c_1$ for some constant $c_1>0$, and achieves  a faster convergence rate of $\mathcal{O}\left(k^{-1}\right)$ in the mean square sense for strongly convex and smooth objective function. To the best of our knowledge, VRA-DGT is the first single-timescale algorithm designed for distributed optimization problems with inexact communication.  Moreover, VRA-DGT achieves convergence rates of $\mathcal{O}(k^{-1})$ in the mean square sense and $\mathcal{O}\left(\frac{\ln(k+1)}{k^b}\right)$, $\forall \in (0.5, 1)$ in the almost sure sense, even under constant stepsizes.
	
	\item[(ii)]For distributed stochastic optimization problems with inexact communication, VRA-DSGT is proposed by combining  VRA-DGT with a hybrid variance-reduced technique \cite{Cutko1019MVR}.  With the integration of the hybrid variance-reduced technique,  VRA-DSGT can further suppress noise from stochastic gradients while maintaining the convergence rate of $\mathcal{O}(k^{-1})$in the mean square sense, which matches the convergence rates of centralized stochastic gradient methods under perfect communication.  Moreover, VRA-DSGT also maintains the convergence rate $\mathcal{O}\left(\frac{\ln(k+1)}{k^b}\right)$, $\forall b\in(0.5,1)$ in the almost sure sense.
Finally, we provide empirical comparisons of our algorithms against existing robust distributed optimization algorithms using real-world data.
\end{itemize}

The rest of this paper is organized as follows.  Section \ref{sec:VRA-DGT} introduces VRA-DGT and VRA-DSGT algorithms for distributed deterministic and stochastic optimization problems respectively.  Section \ref{sec:conv analy} establishes the convergence rate of the proposed algorithms in the mean square and almost sure senses. Finally, Section \ref{sec:num-experi} presents numerical results that validate the effectiveness of the proposed algorithms.

Throughout this paper, we use the following notation. Denote $\mathbb{R}^d$,  $\mathbf{1}$ and $\mathbf{I}$ as the d-dimensional Euclidean space, the vector of ones and the identity matrix respectively. For two matrices $\mathbf{A}$ and $\mathbf{B}$, $\langle \mathbf{A},\mathbf{B}\rangle$  denotes the Frobenius inner product.  $\|\cdot\|_2$  denotes the $l_2$-norm for vector and matrix; 
$\|\cdot\|$ represents the Frobenius norm of matrix. For any two positive sequences $\{a_k\}$ and $\{b_k\}$,  $a_k=\mathcal{O}(b_k)$ if there exists a positive constant $c$ such that $a_k\leq cb_k$.
 Also,  $a_k\asymp b_k$ if $a_k=\mathcal{O}(b_k)$ and $b_k=\mathcal{O}(a_k)$.   %$\rm{Id}$, $[d]$.
Consider $n$ nodes interacting over a graph $\mathcal{G}=(\mathcal{V},\mathcal{E})$, where $\mathcal{V}=\{1,2,\cdots,n\}$ is the set of vertices and $\mathcal{E}\subseteq \mathcal{V}\times \mathcal{V}$ is a collection of
ordered pairs $(i, j)$ such that node $j$ can send information to node $i$. We assume agents interacting over  an undirected graph i.e.,  $(i, j)\in\mathcal{E}$ if and only if  $(j, i)\in\mathcal{E}$. An undirected graph is said to be connected if there exists a path between any two nodes.
%For a graph $\mathcal{G}$, we define a nonnegative weight matrix $\mathbf{W}=[w_{ij}]\in \mathbb{R}^{n\times n}$, where $w_{ij}>0$ $(i, j)\in \mathcal{E}_{\mathbf{W}}$ if and only if , and
Denote $\mathcal{N}_i$ as the collection  of neighbors of agent $i$.  Besides, we assign a weight matrix $\mathbf{W}=[w_{ij}]\in \mathbb{R}^{n\times n}$ to the graph $\mathcal{G}$, where
$w_{ij}>0$ if $(i,j)\in\mathcal{E}$ and $w_{ij}=0$ if  $(i,j)\notin\mathcal{E}$. Moreover, all the proofs of the lemmas are delegated to the Appendix. 

\begin{table}[http]\label{tab:parameters}
	\renewcommand{\tablename}
	\caption{\centering{ Table 1. Symbols in Algorithm 1.}}
	\renewcommand{\arraystretch}{1.2}
	\centering
	\begin{tabular}{l|l}
			\hline
			Symbol & Meaning\\
			\hline
			$x_{i,k}$ & Agent $i$’s local copy of  decision variable at the $k$-th iteration.\cr
			\hline
			$s_{i,k}$ & Agent $i$’s tracker of the cumulative global  gradient $\sum_{t=1}^{k-1}\frac{1}{n}\sum_{j=1}^n\nabla f_j(x_{j,t})$.\cr
			\hline
			$q_{i,k}$ & Local gradient $\nabla f_i(x_{i,k})$ or its  estimator based on stochastic gradients.\cr
			\hline
			$z_{i,k}^s$ &   Estimation of the neighborhood-weighted sum  $\sum_{j \in \mathcal{N}_i} w_{ij} s_{j,k}$.\cr
			\hline
			$z_{i,k}^x$ &  Estimation of the neighborhood-weighted sum  $\sum_{j \in \mathcal{N}_i} w_{ij} x_{j,k}$.\cr
			\hline
			$\tilde{\gamma}$ &  Noise attenuation parameter in the $s_{i,k}$ update.  \cr\hline
			$\gamma_k$ &  Noise attenuation parameter in the $x_{i,k}$ update at $k$-iteration. \cr
			\hline
			$\alpha_k$ & Stepsize at $k$-iteration. \cr
			\hline
			$\lambda_k$ & Parameter of the hybrid variance-reduced technique at $k$-iteration. \cr
			\hline
			$\beta_k$ & VRA mechanism parameter at $k$-iteration. \cr
			\hline
			$\mathbf{W}$&Weight matrix. \cr
			\hline
			$\zeta_{i,k}^s$ & Information-sharing noise in $s_{i,k}$ transmission. \cr
			\hline
			$\zeta_{i,k}^x$ & Information-sharing noise in $x_{i,k}$ transmission. \cr
			\hline
	\end{tabular}
	
\end{table}

\section{VRA-DGT and VRA-DSGT Algorithms}\label{sec:VRA-DGT}
Algorithm \ref{alg:VRA-DSGT} presents the  VRA  based distributed (stochastic) gradient tracking algorithm for distributed (stochastic) optimization problems over the imperfect-communication undirected networks,  where Table \ref{tab:parameters} summarizes the symbols in Algorithm 1.
\renewcommand{\thealgorithm}{1}
\begin{algorithm}[h]
	\caption{\underline{V}ariance-\underline{R}educed \underline{A}ggregation  based \underline{D}istributed (\underline{S}tochastic) \underline{G}radient \underline{T}racking (VRA-DGT and VRA-DSGT)} \label{alg:VRA-DSGT}
	\begin{algorithmic}[1]
		\REQUIRE  initial values $x_{i,1}=z_{i,1}^x\in \mathbb{R}^d$, $s_{i,1}=z_{i,1}^s\in \mathbb{R}^d$, and $q_{i,1}=\nabla f_i(x_{i,1})$ (deterministic case) or $q_{i,1}=\nabla F_i(x_{i,1};\xi_{i,1})$ (stochastic case) for any $i\in\mathcal{V}$; positive factors $\alpha_k,\beta_k,\tilde{\gamma},\gamma_k,\lambda_k\in(0,1]$; nonnegative weight matrix $\mathbf{W}$.
		\FOR {$k=1,2,\cdots$}
		\FOR {$i=1,\cdots,n$ in parallel}
		\STATE Update cumulative global gradient's tracker
		\begin{equation}\label{up:s}
			s_{i,k+1}=\left(1-\tilde{\gamma}\right)s_{i,k}+\tilde{\gamma}\left( z_{i,k}^s+w_{ii}s_{i,k}\right)+q_{i,k}.
		\end{equation} 
		\STATE
		Update decision variable
		\begin{equation}\label{up:x}
			x_{i,k+1}=\left(1-\gamma_k\right)x_{i,k}+\gamma_k\left(z_{i,k}^x+w_{ii}x_{i,k}\right)-\alpha_k \left(s_{i,k+1}-s_{i,k}\right).
		\end{equation} 
		\STATE Update local gradient or its estimator
			\begin{equation}\label{alg:qq}
			q_{i,k+1}=\left\{
			\begin{aligned}
			& \nabla f_i(x_{i,k+1}),\text{ (deterministic case)}\\
			&(1-\lambda_k)\left(q_{i,k}-\nabla F_i(x_{i,k};\xi_{i,k+1})\right)+\nabla F_i(x_{i,k+1};\xi_{i,k+1}).\text{ (stochastic case)}
			\end{aligned}
			\right.
			\end{equation}
		\STATE Agent $i$ receives $\frac{s_{j,k+1}-(1-\beta_k)s_{j,k}}{\beta_k}+\zeta_{j,k+1}^s$, $\frac{x_{j,k+1}-(1-\beta_k)x_{j,k}}{\beta_k}+\zeta_{j,k+1}^x$  from each $j\in \mathcal{N}_{i}$,  and update
			\begin{align}
				&z_{i,k+1}^s=(1-\beta_k)z_{i,k}^s+\beta_k\sum_{j\in \mathcal{N}_{i}} w_{ij}\left(\frac{s_{j,k+1}-(1-\beta_k)s_{j,k}}{\beta_k}+\zeta_{j,k+1}^s\right),\label{alg:z-s}\\
				&z_{i,k+1}^x=(1-\beta_k)z_{i,k}^x+\beta_k\sum_{j\in \mathcal{N}_{i}} w_{ij}\left(\frac{x_{j,k+1}-(1-\beta_k)x_{j,k}}{\beta_k}+\zeta_{j,k+1}^x\right),\label{alg:z-x}
			\end{align}where $\zeta_{j,k+1}^s$ and $\zeta_{j,k+1}^x$ are the information-sharing noise.
		\ENDFOR
		\ENDFOR
	\end{algorithmic}
	\end{algorithm}

In Algorithm \ref{alg:VRA-DSGT}, \textbf{Step 3} updates the estimator $s_{i,k+1}$ of the cumulative global gradient by combining the previous state $s_{i,k}$ with both the estimator $z_{i,k}^s$ of the neighborhood-weighted sum  $\sum_{j \in \mathcal{N}_i} w_{ij} s_{j,k}$ and  the local gradient (estimator) $q_{i,k}$.  Define difference $y_{i,k}\define s_{i,k+1} - s_{i,k}$,  \textbf{Step 3} derives the dynamic-consensus based gradient tracking process
\begin{align}\label{y-s}
		y_{i,k}
		=(1-\tilde{\gamma})y_{i,k-1}+\tilde{\gamma}\sum_{j \in \mathcal{N}_i\cup\{i\}}w_{ij}y_{j,k-1}+\left(q_{i,k}+\tilde{\gamma} e_{i,k}^s\right)-\left(q_{i,k-1}+\tilde{\gamma} e_{i,k-1}^s\right),
\end{align}
where $\tilde{\gamma}$ is a noise attenuation parameter, $e_{i,k}^s:=z_{i,k}^s-\sum_{j\in \mathcal{N}_i} w_{ij} s_{j,k}$ is the estimation error. Indeed,  if  $\tilde{\gamma}=1$,  $q_{i,k}=\nabla f_i(x_{i,k})$ and $e_{i,k}^s=\mathbf{0}$,  (\ref{y-s}) is   exactly the global gradient tracking update of conventional DGT algorithm, which ensures that  $y_{i,k}$ asymptotically approximates the global gradient $\frac{1}{n}\sum_{j=1}^n \nabla f_j(x_{j,k})$. 
 \textbf{Step 4} updates the decision variable $x_{i,k+1}$ by essentially performing a gradient descent step in the direction of $-\frac{1}{n} \sum_{j=1}^n \nabla f_j(x_{j,k})$. Specifically, defining $\bar{x}_k := \frac{1}{n} \sum_{j=1}^n x_{j,k}$, we have
\begin{align}\label{bar-x}
	\bar{x}_{k+1}
	&=\bar{x}_{k}-\alpha_{k}\frac{1}{n}\sum_{j=1}^n\nabla f_j(x_{j,k})+\gamma_k\frac{1}{n}\sum_{j=1}^ne_{j,k}^x+\alpha_{k}\frac{1}{n}\sum_{j=1}^n\left(\tilde{\gamma}e_{j,k}^s+e_{j,k}^q\right),
\end{align}
where $\gamma_k$ is a noise attenuation parameter, $e_{j,k}^q:=q_{j,k}-\nabla f_j(x_{j,k})$, $e_{j,k}^x:=z_{j,k}^x-\sum_{i \in \mathcal{N}_j} w_{ji} x_{i,k}$ along with $e_{j,k}^s$ are the estimation errors. Compared with $\tilde{\gamma}$ in \textbf{Step 3}, $\gamma_k$ in \textbf{Step 4} is time-varying to adapt to  the time-varying stepsize $\alpha_k$. For stochastic gradient case, \textbf{Step 5} updates estimator  $q_{i,k}$  by utilizing stochastic gradients $\nabla F_i(x_{i,k};\xi_{i,k+1})$ and $\nabla F_i(x_{i,k+1};\xi_{i,k+1})$ via hybrid variance-reduced technique \cite{Cutko1019MVR}, which ensures the estimation error $e_{j,k}^q$ in \textbf{Step 4} vanishes asymptotically.
\textbf{Step 6} updates the estimators $z_{i,k+1}^s$ and $z_{i,k+1}^x$ via the VRA mechanism \cite{zhao2023VRA}, which is  a variant of the hybrid variance-reduced technique to suppress information-sharing noise.

 VRA-D(S)GT generalizes both the Robust Push-Pull method \cite{pu2020robust} and the VRA-GT method \cite{zhao2023VRA} in the context of undirected networks.
	Specifically, when $\gamma_k \equiv \tilde{\gamma}$ in (\ref{up:s}), $q_{i,k} = \nabla f_i(x_{i,k})$ in (\ref{up:x}), and $\beta_k \equiv 1$ in both (\ref{alg:z-s}) and (\ref{alg:z-x}), 
	Algorithm \ref{alg:VRA-DSGT} reduces to	the Robust Push-Pull method; when $q_{i,k} = \nabla f_i(x_{i,k})$ in (\ref{alg:qq}) and $\beta_k \equiv 1$ in (\ref{alg:z-x}), Algorithm \ref{alg:VRA-DSGT} reduces to the
	VRA-GT method. Moreover, VRA-D(S)GT is also an extension of the GT-HSGD method \cite{xin21Hybrid} which incorporates the hybrid variance-reduced technique into the DGT algorithm and focuses on the perfect communication case.

For the sake of analysis, define variables
\begin{align}
	&\mathbf{x}_k:= [x_{1,k},x_{2,k},\cdots,x_{n,k}]^\intercal,~\mathbf{y}_k:= [s_{1,k+1}-s_{1,k},s_{2,k+1}-s_{2,k},\cdots,s_{n,k+1}-s_{n,k}]^\intercal,\label{x}\\
	&\mathbf{z}_k^s:= [z_{1,k}^s,z_{2,k}^s,\cdots,z_{n,k}^s]^\intercal,~\mathbf{z}_k^x:= [z_{1,k}^x,z_{2,k}^x,\cdots,z_{n,k}^x]^\intercal,\label{z}\\
	&\mathbf{q}_k= [q_{1,k},q_{2,k},\cdots,q_{n,k}]^\intercal,~\mathbf{g}_k:=\left[\nabla f_1(x_{1,k}),\nabla f_2(x_{2,k}),\cdots,\nabla f_n(x_{n,k})\right]^\intercal,
	\notag\\
	&\tilde{\mathbf{g}}_k= [\nabla F_1(x_{1,k};\xi_{1,k+1}),\nabla F_2(x_{2,k};\xi_{2,k+1}),\cdots,\nabla F_n(x_{n,k};\xi_{n,k+1})]^\intercal,\notag\\
	&\tilde{\mathbf{g}}_{k+1}^{'}= [\nabla F_1(x_{1,k+1};\xi_{1,k+1}),\nabla F_2(x_{2,k+1};\xi_{2,k+1}),\cdots,\nabla F_n(x_{n,k};\xi_{n,k+1})]^\intercal,\notag
\end{align}
weight matrices
\begin{align*}
	&\mathbf{W}_{\tilde{\gamma}}:=(1-\tilde{\gamma})\mathbf{I}+\tilde{\gamma}\mathbf{W},~\mathbf{W}_{\gamma_k}:=(1-\gamma_k)\mathbf{I}+\gamma_k\mathbf{W},%~\mathbf{W}_{d}=\diag\left(w_{11},w_{22},\dots,w_{nn}\right),
	~\tilde{\mathbf{W}}:=\mathbf{W}-\diag\left(w_{11},\dots,w_{nn}\right)\label{w}
\end{align*}
and estimation errors
\begin{equation}
	\mathbf{e}_k^x:=\mathbf{z}^x_k-\tilde{\mathbf{W}}\mathbf{x}_k,~\mathbf{e}_k^s:=\mathbf{z}^s_k-\tilde{\mathbf{W}}\mathbf{s}_k.\label{e}	
\end{equation}
Then, VRA-D(S)GT has the matrix form of
\begin{align}
	&\mathbf{x}_{k+1}=\mathbf{W}_{\gamma_k}\mathbf{x}_k+\gamma_k\mathbf{e}_k^x-\alpha_k\mathbf{y}_k\label{alg:x}\\
	&\mathbf{q}_{k+1}=\left\{
	\begin{aligned}
	&\mathbf{g}_{k+1}\text{ (deterministic case)}\\
	&(1-\lambda_k)\left(\mathbf{q}_k-\tilde{\mathbf{g}}_k\right)+\tilde{\mathbf{g}}_{k+1}^{'}\text{ (stochastic case)}
	\end{aligned}
	\right.\label{alg:q}\\
	&\mathbf{y}_{k+1}=\mathbf{W}_{\tilde{\gamma}}\mathbf{y}_{k}+\left(\mathbf{q}_{k+1}+\tilde{\gamma}\mathbf{e}_{k+1}^s\right)-\left(\mathbf{q}_{k}+\tilde{\gamma}\mathbf{e}_{k}^s\right).\label{alg:y}
\end{align}

\section{Convergence Analysis}\label{sec:conv analy}

In this section, we provide the  convergence rates of VRA-D(S)GT algorithm.  We first collect all the assumptions  and provide some technical results that will be used throughout the paper.

\begin{ass}[\textbf{objective function}]\label{ass:function}
	The local functions $f_i(x),~i\in\mathcal{V}$ are $L$-smooth and the global function $f(x)$ is $\mu$-strongly convex, i.e., for any $x,y \in \mathbb{R}^d$,
	\begin{align}\label{g-l}
		&\|\nabla f_i(x)-\nabla f_i(y)\|\le  L\|x-y\|
	\end{align}
	and
	\begin{align*}
		&f(y)\ge f(x)+\left\langle \nabla f(x),y-x\right\rangle+\frac{\mu}{2}\|x-y\|^2.
	\end{align*}
\end{ass}
The L-smoothness of $f_i$ can be immediately translated to the  L-smoothness of $f(x)$. The strong convexity of $f(x)$ can guarantee that problem (\ref{model}) has a unique optimal solution $x^*$.

\begin{ass}  [\textbf{networks}]\label{ass:matrix}
	The graph $\mathcal{G}$ corresponding to the network of agents is undirected and
	connected. The weight matrix $\mathbf{W}$ assigned to $\mathcal{G}$ is doubly stochastic, i.e., $\mathbf{W}\mathbf{1}=\mathbf{1}$ and $\mathbf{1}^\intercal\mathbf{W}=\mathbf{1}^\intercal$. In addition, $w_{ii}>0$ for some $i\in\mathcal{V}$.
\end{ass}
Assumption \ref{ass:matrix} implies  that there exists a constant $\eta_w\in(0,1]$ such that
\begin{equation}\label{eta}
	\left\|(1-\gamma)\mathbf{I}+\gamma\mathbf{W}-\frac{\mathbf{1}\mathbf{1}^\intercal}{n}\right\|_2\le 1-\gamma \eta_w
\end{equation}
for any $\gamma\in [0,1]$ \cite[Section II-B]{qu2017harnessing}.

We work with a rich enough probability space $\left(\Omega,\mathcal{F},\mathbb{P}\right)$ and define two natural filtration as
\begin{align}
& \mathcal{F}^x_k\define \sigma\left(\left\{\xi_{j,t},\zeta_{j,t}^x,\zeta_{j,t}^s:2\le t\le k-1,~ j\in\mathcal{V} \right\}\right)~(\forall k\ge 3),~\mathcal{F}^x_1=\mathcal{F}^x_2=\left\{\Omega,\emptyset\right\}
\end{align} 
and 
\begin{align}
	& \mathcal{F}^s_k\define \sigma\left(\left\{\xi_{j,t},\zeta_{j,t}^x,\zeta_{j,t}^s:2\le t\le k,~ j\in\mathcal{V} \right\} \backslash \left\{\zeta_{j,k}^x:j\in\mathcal{V} \right\}\right)~(\forall k\ge 2),~\mathcal{F}^s_1=\left\{\Omega,\emptyset\right\},
\end{align}
where $\emptyset$ is the empty set. Obviously, $\mathcal{F}_1^x=\mathcal{F}_1^s\subseteq \mathcal{F}_2^x\subseteq\mathcal{F}_2^s\subseteq \mathcal{F}_3^x\subseteq\mathcal{F}_3^s\cdots$.

\begin{ass}[\textbf{information-sharing noise}]\label{ass:noise}
	For any  $k\ge 1$, $i, j\in\mathcal{V}$ and $i\ne j$, $$\mathbb{E}\left[\zeta^s_{i,k}\big|\mathcal{F}_k^x\right]=\mathbf{0},~\mathbb{E}\left[\zeta_{i,k}^x\big|\mathcal{F}_k^s\right]=\mathbf{0},~\mathbb{E}\left[\langle\zeta^s_{i,k}, \zeta^s_{j,k}\rangle\big|\mathcal{F}_k^x\right]=0.$$ 
	Moreover, there exist constants $\sigma_\zeta>0$, $p_1\ge 2$ such that $ \mathbb{E}\left[\|\zeta^s_{i,k}\|_2^{p_1}\big|\mathcal{F}_k^x\right]\le\sigma_\zeta^{p_1}$ and  $ \mathbb{E}\left[\|\zeta^x_{i,k}\|_2^{p_1}\big|\mathcal{F}_k^s\right]\le\sigma_\zeta^{p_1}$.
\end{ass}

\begin{ass}[\textbf{stochastic gradient}]\label{ass:stochastic gradient}
	For any $i\in\mathcal{V}$, there exist constants $\sigma_\xi>0$, $p_2\ge 2$ and random variable $L_i(\xi_i)$ such that
	\begin{align}\label{sg-var}
		&\mathbb{E}\left[\nabla F_i(x;\xi_i)\right]=\nabla f_i(x),~\mathbb{E}\left[\left\|\nabla F_i(x;\xi_i)-\nabla f_i(x)\right\|_2^{p_2}\right]\le  \sigma_\xi^{p_2}
	\end{align}
	and
	\begin{align}\label{sg-l}
		\left\|\nabla F_i(x;\xi_i)-\nabla F_i(y;\xi_i)\right\|_2 \le L_i(\xi_i) \|x-y\|_2, \forall x,y\in\mathbb{R}^d.
	\end{align}
	In addition,  $L^2\define\max_{i\in\mathcal{V}} \mathbb{E}\left[ L_i(\xi_i)^2\right]<\infty$\footnote{Note that (\ref{sg-l}) implies (\ref{g-l}) by Jensen's inequality, that is $\|\nabla f_i(x)-\nabla f_i(y)\|\le  \sqrt{\mathbb{E}\left[ L_i(\xi_i)^2\right]} \|x-y\|$. Hence, despite a slight confusion in notation, we simply denote $L\define\sqrt{\mathbb{E}\left[ L_i(\xi_i)^2\right]}$ here.}.
	
\end{ass}

When $p_1=p_2=2$, Assumptions \ref{ass:noise} and \ref{ass:stochastic gradient} represent standard assumptions for distributed stochastic optimization over the imperfect-communication networks. Specifically, Assumption \ref{ass:noise} holds in various scenarios, including noisy communication channels \cite{Dasarathan2015robust,Kar2009imc} and communication quantization \cite{Doan2021Quanti,YUAN2012DDA}. Assumption \ref{ass:stochastic gradient} implies the Lipschitz smoothness of each $f_i$ and ensures that the stochastic gradient is unbiased and has bounded variance. On the other hand, the case where $p_1>2$ and $p_2>2$ is required to establish the almost sure convergence rate of VRA-D(S)GT.

The following propositions provide some  technical results on the VRA mechanism.

\begin{prop}\label{prop:vra}
Suppose Assumption \ref{ass:noise} (with $p_1=2$) holds. Then, there exist constants $c_x$ and $c_s$ such that
\begin{equation}\label{es-b1}
	\mathbb{E}\left[\left\|\mathbf{e}^s_{k}\right\|^2\right]\le c_s\beta_k,~\mathbb{E}\left[\left\|\mathbf{e}^x_{k}\right\|^2\right]\le c_x\beta_k,
\end{equation}	 
when $\beta_k\asymp \frac{1}{k^b}, b\in (1/2,1), $ or $\beta_k=\frac{a_2}{k+a_1}, a_1\ge 0, a_2\in[1,a_1+1)$.
\end{prop}
\begin{proof} 
	 According to the definition of $\mathbf{e}^s_{k}$ in (\ref{e}) and recursion of $z_{i,k+1}^s$ in (\ref{alg:z-s}),
	\begin{align}\label{e-s}
		\mathbf{e}^s_{k+1}
		&=(1-\beta_k)\left(\mathbf{z}_{k}^s+\tilde{\mathbf{W}}\mathbf{s}_{k+1}-\tilde{\mathbf{W}}\mathbf{s}_{k}\right)+\beta_k \tilde{\mathbf{W}}\left(\mathbf{s}_{k+1}+\zeta_{k+1}^s\right)-\tilde{\mathbf{W}} \mathbf{s}_{k+1}\notag\\
		&=(1-\beta_k)\mathbf{e}^s_{k}+\beta_k \tilde{\mathbf{W}}\zeta_{k+1}^s,
	\end{align}
	where $\zeta_{k+1}^s:=\left[\zeta_{1,k+1}^s,\zeta_{2,k+1}^s,\cdots,\zeta_{n,k+1}^s\right]^\intercal$. Without loss of generality, we let $\mathbf{e}^s_1=\mathbf{0}$.
	By induction, we have the following decomposition
	\begin{equation}\label{e-recur}
		\mathbf{e}^s_{k+1}=\Pi_{t=1}^k(1-\beta_t)\mathbf{e}^s_1+\sum_{t=1}^k(\Pi_{l=t+1}^k(1-\beta_l))\beta_t\tilde{\mathbf{W}}\zeta_{t+1}^s=\sum_{t=1}^k(\Pi_{l=t+1}^k(1-\beta_l))\beta_t\tilde{\mathbf{W}}\zeta_{t+1}^s,
	\end{equation}
	where $\Pi_{l=t+1}^k(1-\beta_l)=1$ for $t=k$. Then, when $\beta_{k}=\frac{1}{k+a_1}$,
	\begin{align*}
		\mathbb{E}\left[\|\mathbf{e}^s_{k+1}\|^2\right]&=	\mathbb{E}\left[\left\|\frac{1}{k+a_1}\sum_{t= 1}^k\tilde{\mathbf{W}}\zeta_{t+1}^s\right\|^2\right]\\
		%&=\frac{1}{(k+a_1)^2}\sum_{t_1= 1}^k\sum_{t_2= 1}^k\mathbb{E}\left[\text{Tr}\left((\tilde{\mathbf{W}}\zeta_{t_1+1}^s)^\intercal\tilde{\mathbf{W}}\zeta_{t_2+1}^s\right)\right]\\
		&=\frac{1}{(k+a_1)^2}\sum_{t_1= 1}^k\sum_{t_2= 1}^k\mathbb{E}\left[\mathbb{E}\left[\text{Tr}\left((\tilde{\mathbf{W}}\zeta_{t_1+1}^s)^\intercal\tilde{\mathbf{W}}\zeta_{t_2+1}^s\right)\big|\mathcal{F}_{\max \{t_1,t_2\}+1}^x\right]\right]\\
		&=\frac{1}{(k+a_1)^2}\sum_{t=1}^k\mathbb{E}\left[\left\|\tilde{\mathbf{W}}\zeta_{t+1}^s\right\|^2\right]\\
		%&\le \frac{n\sigma_\zeta^2}{k+a_1}\\
		&\le c_s\beta_{k+1},
	\end{align*}
	where $c_s=2n\sigma_\zeta^2$, $\text{Tr}(\mathbf{A})$ represents the trace of matrix $\mathbf{A}$, the third equality follows from the condition $\mathbb{E}\left[\zeta_{i,k}^s\big|\mathcal{F}_k^x\right]=\mathbf{0}$ for any $k\ge 1$, the  inequality follows from the condition $ \mathbb{E}\left[\|\zeta^s_{i,k}\|_2^2\big|\mathcal{F}_k^x\right]\le\sigma_\zeta^2$. 
	
	For the cases of $\beta_k\asymp \frac{1}{k^b}, b\in (1/2,1), $ or $\beta_k=\frac{a_2}{k+a_1}, a_1\ge 0, a_2\in(1,a_1+1)$, see  {\cite[Theorem 1]{zhao2023VRA}}.
	
	By a similarly analysis, we can show that there exists constant $c_x$ such that
$$
	\mathbb{E}\left[\left\|\mathbf{e}^x_{k}\right\|^2\right]\le c_x\beta_k.$$
The proof is complete. 
\end{proof}

\begin{prop}\label{prop:vra-1}
		Suppose Assumption \ref{ass:noise} (with $p_1>2$) holds. Then, there exist positive finite random variables $\tilde{c}_x$ and $\tilde{c}_s$ such that 
		$$
		\left\|\mathbf{e}^s_{k}\right\|^2\le \tilde{c}_s\beta_k\ln (k+1),~\left\|\mathbf{e}^x_{k}\right\|^2\le \tilde{c}_x\beta_k\ln (k+1)~~~\text{a.s.}
		$$
		when $\beta_k\asymp\frac{1}{k^b}$, $b\in\left(p_1/(2p_1-2),1\right)$.
\end{prop}
\begin{proof}
		  By  (\ref{e-recur}), it is equivalent to show 
	\begin{equation}\label{sum-as}
		\sum_{t=1}^k(\Pi_{l=t+1}^k(1-\beta_l))\beta_t\tilde{\mathbf{W}}\zeta_{t+1}^s=\mathcal{O}\left(\beta_k\ln (k+1)\right)~~~\text{a.s.}
	\end{equation}
	In what follows, we show (\ref{sum-as}) by Lemma \ref{lem:as} in Appendix B.	Denote
	\begin{equation*}
		\Gamma=\mathbf{I},~R_k=\tilde{\mathbf{W}},~\tilde{\zeta}_k=\zeta_{k}^s,~\mathcal{F}_k=\mathcal{F}_{k+1}^x,~v_k=1.
	\end{equation*}
	Then (\ref{sum-as}) can be rewritten as
	\begin{equation*}
		M_{k+1}=\sum_{t=1}^kB_{k,t}\beta_{t}R_t\tilde{\zeta}_{t+1},
	\end{equation*}
	which is in the form of (\ref{M-def}) in Lemma \ref{lem:as}. Next, we verify the conditions (C1)-(C4) of Lemma \ref{lem:as}. 
	
	According to the definitions of $\beta_{k}$ and $\tilde{\mathbf{W}}$, conditions (C2)-(C4) of Lemma \ref{lem:as} hold. 
	
	Noting that $\{\zeta_{k}^s\}$ is a martingale difference sequence, it remains to verify the following inequalities of the condition (C1) of Lemma \ref{lem:as}: 
	\begin{align*}
		&\mathbb{E}\left[\|\tilde{\zeta}_{k+1}\|^2 \big|\mathcal{F}_k\right]\le C+R_{2,k}~a.s.,\\
		&\sum_{t= 1}^\infty \beta_t\mathbb{E}\left[\|\tilde{\zeta}_{t+1}\|^21_{\left\{\|\tilde{\zeta}_{t+1}\|^2\ge \frac{1}{\beta_t\ln (t+1)}\right\}}\big|\mathcal{F}_k\right]< \infty,
	\end{align*}
	where  $C$ is some positive constant and  $\{R_{2,k}\}$ is a sequence of positive random variables converging almost surely to 0.
	Assumption \ref{ass:noise} implies
	\begin{align*}
		\mathbb{E}\left[\|\tilde{\zeta}_{k+1}\|^2 \big|\mathcal{F}_k\right]=\mathbb{E}\left[\|\zeta_{k+1}^s\|^2 \big|\mathcal{F}_{k+1}^x\right]\le n\sigma_\zeta^2.
	\end{align*}
	Noting that $	\mathbb{E}\left[\|\zeta^s_{i,k+1}\|^{p_1} \big|\mathcal{F}_{k+1}^x\right]\le \sigma_\zeta^{p_1}$ where $p_1>2$ is defined in Assumption \ref{ass:noise}, then by   H\"{o}lder inequality and  Markov inequality, we obtain
	\begin{align*}
		\mathbb{E}\left[\|\tilde{\zeta}_{t+1}\|^21_{\left\{\|\tilde{\zeta}_{t+1}\|^2\ge \frac{1}{\beta_t\ln (t+1)}\right\}}\big|\mathcal{F}_t\right]&=\mathbb{E}\left[\|\zeta_{t+1}^s\|^21_{\left\{\|\zeta^s_{t+1}\|^2\ge \frac{1}{\beta_t\ln (t+1)}\right\}}\big|\mathcal{F}_{t+1}^x\right]\\
		&\le \left(\mathbb{E}\left[\|\zeta_{t+1}^s\|^{p_1}\big|\mathcal{F}_{t+1}^x\right]\right)^{2/p_1}\left(\mathbb{E}\left[1^{q_1}_{\left\{\|\zeta^s_{t+1}\|^2\ge \frac{1}{\beta_t\ln (t+1)}\right\}}\big|\mathcal{F}_{t+1}^x\right]\right)^{1/q_1}\\
		&= \left(\mathbb{E}\left[\|\zeta_{t+1}^s\|^{p_1}\big|\mathcal{F}_{t+1}^x\right]\right)^{2/p_1}\left(\mathbb{P}\left[\|\zeta^s_{t+1}\|^2\ge \frac{1}{\beta_t\ln (t+1)}\big|\mathcal{F}_{t+1}^x\right]\right)^{1/q_1}\\
		&\le \left(\mathbb{E}\left[\|\zeta_{t+1}^s\|^{p_1}\big|\mathcal{F}_{t+1}^x\right]\right)^{2/p_1}\left(\mathbb{E}\left[\|\zeta_{t+1}^s\|^2\big|\mathcal{F}_{t+1}^x\right]\beta_t\ln (t+1)\right)^{1/q_1}\\
		&\le \sigma_\zeta^{2+2/q_1}\left(\beta_t\ln (t+1)\right)^{1/q_1},
	\end{align*}
	where $q_1$ is the constant satisfying $\frac{2}{p_1}+\frac{1}{q_1}=1$.
	Subsequently,
	\begin{align*}
		\sum_{t= 1}^\infty \beta_t\mathbb{E}\left[\|\tilde{\zeta}_{t+1}\|^21_{\left\{\|\tilde{\zeta}_{t+1}\|^2\ge \frac{1}{\beta_t\ln (t+1)}\right\}}\big|\mathcal{F}_t\right]\le \sum_{t= 1}^\infty \beta_t^{1+1/q_1} \sigma_\zeta^{2+2/q_1}\left(\ln (t+1)\right)^{1/q_1}< \infty,
	\end{align*}
	where the last inequality follows from the fact $b>\frac{p_1}{2(p_1-1)}$.
	Then by Lemma \ref{lem:as}, $$\left\|\mathbf{e}^s_{k+1}\right\|^2=\left\|\sum_{t=1}^k(\Pi_{l=t+1}^k(1-\beta_l))\beta_t\tilde{\mathbf{W}}\zeta_{t}^s\right\|^2=\mathcal{O}\left(\beta_k\ln (k+1)\right)$$ almost surely, which combining with the fact $\lim_{k\rightarrow\infty}\frac{\beta_k}{\beta_{k+1}}<\infty$ implies $\left\|\mathbf{e}^s_{k}\right\|^2\le \tilde{c}_s\beta_k\ln (k+1)$ for some positive finite random variable $\tilde{c}_s$.
	
	By a similar analysis, 	we can show that there exists a positive finite random variable $\tilde{c}_x$ such that 
	$$\left\|\mathbf{e}^x_{k}\right\|^2\le \tilde{c}_x\beta_k\ln (k+1)~~~\text{a.s.}
	$$
	The proof is complete. 
\end{proof}

The following lemmas   provide the upper bounds of  hybrid variance-reduced error, gradient tracking error, agreement error and  optimality gap of the algorithm.

\begin{lem}[\textbf{hybrid variance-reduced error}]\label{lem:vr}Let
	\begin{equation}\label{eq}
		\mathbf{e}_{k+1}^q=\mathbf{q}_{k+1}-\mathbf{g}_{k+1}.
	\end{equation}
	Then under Assumptions \ref{ass:function}-\ref{ass:stochastic gradient},
	\begin{align*}
		\mathbb{E}\left[\left\|\mathbf{e}_{k+1}^q\right\|^2\big|\mathcal{F}_k^s\right]
	&\le (1-\lambda_k)^2\left\|\mathbf{e}_{k}^q\right\|^2+30L_\xi^2(1-\lambda_k)^2\left(\gamma_k^2\left\|\mathbf{W}-\mathbf{I}\right\|^2\left\|\mathbf{x}_k-\bar{\mathbf{x}}_k\right\|^2+\gamma_k^2\mathbb{E}\left[\left\|\mathbf{e}_k^x\right\|^2\big|\mathcal{F}_k^s\right]\right.\notag\\
	&\left.\quad+\alpha_k^2\left\|\mathbf{y}_k-\bar{\mathbf{y}}_k\right\|^2+\alpha_k^2n\left\|\bar{q}_k\right\|_2^2+5\alpha_k^2\tilde{\gamma}^2\left\|\mathbf{e}_k^s\right\|^2\right)+3\lambda_k^2n\sigma_\xi^2,
	\end{align*}
	where $\bar{q}_{k}:=\frac{1}{n} \sum_{j=1}^nq_{j,k}$.
\end{lem}
\begin{proof}
See Appendix A.1.
\end{proof}

\begin{lem}[\textbf{gradient tracking error}]\label{lem:y-1}
	Under Assumptions \ref{ass:function}-\ref{ass:stochastic gradient},
\begin{align*}
	&\mathbb{E}\left[\left\|\mathbf{y}_{k+1}-\bar{\mathbf{y}}_{k+1}\right\|^2\big|\mathcal{F}_k^s\right]\\
&\le \left(1-\tilde{\gamma} \eta_w+\frac{40L^2\left(1+2\lambda_k^2\right)}{\tilde{\gamma} \eta_w}\alpha_k^2\right)\left\|\mathbf{y}_{k}-\bar{\mathbf{y}}_{k}\right\|^2+\frac{40L^2\left(1+2\lambda_k^2\right)}{\tilde{\gamma} \eta_w}\left(\gamma_k^2\mathbb{E}\left[\left\|\mathbf{e}_k^x\right\|^2\big|\mathcal{F}_k^s\right]\right.\notag\\
&\quad\left.+\gamma_k^2\left\|\mathbf{W}-\mathbf{I}\right\|^2\left\|\mathbf{x}_k-\bar{\mathbf{x}}_k\right\|^2+n\alpha_k^2\left\|\bar{q}_k\right\|_2^2\right)+\frac{8\lambda_k^2}{\tilde{\gamma} \eta_w}\left(\left\|\mathbf{e}_k^q\right\|^2+n\sigma_\xi^2\right)\notag\\
&\quad+\frac{4\tilde{\gamma}\left(\beta_k^2+10L^2\left(1+2\lambda_k^2\right)\alpha_k^2\right)}{\eta_w}\left\|\mathbf{e}_{k}^s\right\|^2+\frac{4\tilde{\gamma}\beta_k^2}{ \eta_w}n\sigma_\zeta^2,
\end{align*}
	where  $\eta_w$ is defined in (\ref{eta}) and	 $\bar{q}_{k}:=\frac{1}{n} \sum_{j=1}^nq_{j,k}$.
\end{lem}
\begin{proof}
	See Appendix A.2.
\end{proof}

\begin{lem}[\textbf{agreement error}]\label{lem:x-1}
	Under Assumptions \ref{ass:function}-\ref{ass:stochastic gradient},	
		\begin{align*}
		\mathbb{E}\left[\left\|\mathbf{x}_{k+1}-\bar{\mathbf{x}}_{k+1}\right\|^2\big|\mathcal{F}_k^s\right]
		&\le (1-\gamma_k \eta_w)\left\|\mathbf{x}_{k}-\bar{\mathbf{x}}_{k}\right\|^2+\frac{2\alpha_k^2}{\gamma_k \eta_w}\left\|\mathbf{y}_{k}-\bar{\mathbf{y}}_{k}\right\|^2+\frac{2\gamma_k}{ \eta_w}\mathbb{E}\left[\left\|\mathbf{e}_{k}^x\right\|^2\big|\mathcal{F}_k^s\right],
	\end{align*}
	where  $\eta_w$ is defined in (\ref{eta}).
\end{lem}
\begin{proof}
	See Appendix A.3.
\end{proof}

\begin{lem}[\textbf{optimality gap}]\label{lem:xstar-1}
	Under Assumptions \ref{ass:function}-\ref{ass:stochastic gradient},
	\begin{align*}
		\mathbb{E}\left[f(\bar{x}_{k+1})-f^*\big|\mathcal{F}_k^s\right]	
		&\le (1-\mu\alpha_k)(f(\bar{x}_k)-f^*)-\alpha_k\|\bar{q}_{k}\|_2^2+\frac{8\alpha_k}{n}\left\|\mathbf{e}^q_k\right\|^2+\frac{8\alpha_kL^2}{n}\left\|\mathbf{x}_k-\bar{\mathbf{x}}_k\right\|^2\\
		&\quad+\frac{10\alpha_k\tilde{\gamma}^2}{n}\|\mathbf{e}^s_{k}\|^2+\frac{10\gamma_k^2}{\alpha_kn}\mathbb{E}\left[\|\mathbf{e}^x_{k}\|^2\big|\mathcal{F}_k^s\right],
	\end{align*}
	where $f^*=\min_{x\in\mathbb{R}^d} f(x)$, 	 $\bar{q}_{k}:=\frac{1}{n} \sum_{j=1}^nq_{j,k}$.
\end{lem}
\begin{proof}
	See Appendix A.4.
\end{proof}

%The following lemma  provides the upper bound of the distance from averaged iterates to the optimal solution.

\subsection{Convergence Rate in the Mean Square Sense}

In this subsection, we provide the convergence rate of VRA-D(S)GT in the mean square sense.
	To measure the convergence rates of algorithms, we define  Lyapunov functions 
	\begin{equation*}\label{L-func-1}
		\mathcal{L}^S_{k}:=\left\|\mathbf{y}_{k}-\bar{\mathbf{y}}_{k}\right\|^2+\left\|\mathbf{x}_{k}-\bar{\mathbf{x}}_{k}\right\|^2+(f(\bar{x}_k)-f^*)+\left\|\mathbf{e}_k^q\right\|^2
	\end{equation*} 
	for VRA-DSGT, and
	\begin{equation}\label{L-fun-2}
		\mathcal{L}^D_{k}:=\left\|\mathbf{y}_{k}-\bar{\mathbf{y}}_{k}\right\|^2+\left\|\mathbf{x}_{k}-\bar{\mathbf{x}}_{k}\right\|^2+(f(\bar{x}_k)-f^*)
	\end{equation}
	for VRA-DGT.

\begin{thm}\label{thm:x-opt-mse}
	Suppose Assumptions \ref{ass:function}, \ref{ass:matrix}, \ref{ass:noise} (with $p_1=2$) and \ref{ass:stochastic gradient} (with $p_2=2$) hold. Suppose also that $\alpha_k$ is nonincreasing and diminishes to zero, 
	\begin{equation}\label{para-0}
		\gamma_k = c_1 \alpha_k, \quad \lambda_k = c_2 \alpha_k, \quad  \beta_k \le c_3 \alpha_k,
	\end{equation}
	where the parameters $\tilde{\gamma}$, $c_1$, $c_2$, and $c_3$ satisfy
	\begin{align}\label{para}
		\tilde{\gamma} \in (0,1], ~\delta^S\define\min\left\{\mu,\frac{\tilde{\gamma} \eta_w}{\alpha_1}-\frac{2}{c_1 \eta_w},c_1 \eta_w-\frac{8L^2}{n},2c_2-\frac{8}{n}\right\}>0,~	c_2\le \frac{1}{\sqrt{2}\alpha_1},~c_3>0,
	\end{align}
	and	$\eta_w$ is defined in (\ref{eta}).
Then for VRA-DSGT,
	\begin{itemize}
		\item[(i)]
			$\mathbb{E}\left[\mathcal{L}^S_{k+1}\right]=\mathcal{O}\left(\frac{1}{k^{b}}\right)$ 
			when $\alpha_k\asymp\frac{1}{k^{b}}$
			with $b\in(1/2,1)$;
		\item[(ii)]	
		$\mathbb{E}\left[\mathcal{L}^S_{k+1}\right]=\mathcal{O}\left(\frac{1}{k}\right)$ 
		when  
		$\alpha_k=\frac{a_2}{k+a_1}$, $\beta_k=\frac{a_3}{k+a_1}$ with $a_1>0$, $a_2> \frac{2}{\delta^S}
		$  and $a_3\in[1,a_1+1)$.
	\end{itemize}
\end{thm}

\begin{proof}
\textbf{Part (i)}.  By Lemmas \ref{lem:vr}-\ref{lem:xstar-1} and the fact $1+2\lambda_k^2\le 2$, we have
{\small\begin{align}
		\mathbb{E}\left[\mathcal{L}^S_{k+1}\big|\mathcal{F}_k^s\right]
		\le&(1-\mu\alpha_k)(f(\bar{x}_k)-f^*)+\left[1-\tilde{\gamma} \eta_w+\frac{2\alpha_k^2}{\gamma_k \eta_w}+\left(\frac{80L^2}{\tilde{\gamma} \eta_w}+30L^2\right)\alpha_k^2\right]\left\|\mathbf{y}_{k}-\bar{\mathbf{y}}_{k}\right\|^2\notag\\
		&+\left[1-\gamma_k \eta_w+\frac{8\alpha_kL^2}{n}+\left(\frac{80L^2}{\tilde{\gamma} \eta_w}+30L^2\right)\left\|\mathbf{W}-\mathbf{I}\right\|^2\gamma_k^2\right]\left\|\mathbf{x}_k-\bar{\mathbf{x}}_k\right\|^2\notag\\
		&+\left[(1-\lambda_k)^2+\frac{8\lambda_k^2}{\tilde{\gamma} \eta_w}+\frac{8\alpha_k}{n}\right]\left\|\mathbf{e}_k^q\right\|^2\notag\\
		&+\left[\left(\frac{80L^2}{\tilde{\gamma} \eta_w}+30L^2\right)\gamma_k^2+\frac{2\gamma_k}{ \eta_w}+\frac{10\gamma_k^2}{\alpha_kn}\right]\mathbb{E}\left[\left\|\mathbf{e}_k^x\right\|^2\big|\mathcal{F}_k^s\right]\notag\\
		&+\left[150L^2\alpha_k^2\tilde{\gamma}^2+\frac{4\tilde{\gamma}\beta_k^2+80\tilde{\gamma}L^2\alpha_k^2}{\eta_w}+\frac{10\alpha_k\tilde{\gamma}^2}{n}\right]\left\|\mathbf{e}_k^s\right\|^2\notag\\
		& -\left[\alpha_k-\left(\frac{80L^2}{\tilde{\gamma} \eta_w}+30L^2\right)n\alpha_k^2\right]\|\bar{q}_{k}\|_2^2 +\frac{8\lambda_k^2}{\tilde{\gamma} \eta_w}\left(n\sigma_\xi^2\right)+\frac{4\tilde{\gamma}\beta_k^2}{ \eta_w}n\sigma_\zeta^2+3\lambda_k^2n\sigma_\xi^2.\label{L-bound-2}
\end{align}}
Taking the settings  (\ref{para-0}) and (\ref{para}) into (\ref{L-bound-2}),  we have
\begin{align}
\mathbb{E}\left[\mathcal{L}^S_{k+1}\big|\mathcal{F}_k^s\right]
\le&\left(1-\delta^S\alpha_k+\tilde{\delta}^S\alpha_k^2\right)\mathcal{L}^S_{k}-\tau^S_k\|\bar{q}_{k}\|_2^2+\mathbb{E}\left[\tilde{\tau}^S_k\big|\mathcal{F}_k^s\right],\label{L-bound-3}
\end{align}
where $\delta^S$ is defined in (\ref{para}),
\begin{align*}
\tilde{\delta}^S:=&\max\left\{\left(\frac{80L^2}{\tilde{\gamma} \eta_w}+30L^2\right),\left(\frac{80L^2}{\tilde{\gamma} \eta_w}+30L^2\right)\left\|\mathbf{W}-\mathbf{I}\right\|^2c_1^2,\left(1+\frac{8}{\tilde{\gamma} \eta_w}\right)c_2^2\right\}\\
\tau^S_k:=&\alpha_k-\left(\frac{80L^2}{\tilde{\gamma} \eta_w}+30L^2\right)n\alpha_k^2,\\
\tilde{\tau}^S_k:=&\left[\left(\frac{80L^2}{\tilde{\gamma} \eta_w}+30L^2\right)\gamma_k^2+\frac{2\gamma_k}{ \eta_w}+\frac{10\gamma_k^2}{\alpha_kn}\right]\left\|\mathbf{e}_k^x\right\|^2 +\frac{8\lambda_k^2}{\tilde{\gamma} \eta_w}\left(n\sigma_\xi^2\right)+\frac{4\tilde{\gamma}n\sigma_\zeta^2\beta_k^2}{ \eta_w}+3n\sigma_\xi^2\lambda_k^2\notag\\
&+\left[150L^2\alpha_k^2\tilde{\gamma}^2+\frac{4\tilde{\gamma}\left(\beta_k^2+20L^2\alpha_k^2\right)}{\eta_w}+\frac{10\alpha_k\tilde{\gamma}^2}{n}\right]\left\|\mathbf{e}_k^s\right\|^2.
\end{align*}
Noting that  $\alpha_k$    diminishes to zero, 
% we have
%\begin{equation*}
%\alpha_k\le \frac{\delta^S}{2\tilde{\delta}^S},~\tau^S_k\ge 0
%\end{equation*}when $k$ is large enough. 
there exists a large enough integer $k_0$  such that 
\begin{align}\label{L-bound-4}
		\mathbb{E}\left[\mathcal{L}^S_{k+1}\big|\mathcal{F}_k^s\right]
		\le&\left(1-\frac{\delta^S}{2}\alpha_k\right)\mathcal{L}^S_{k}+\mathbb{E}\left[\tilde{\tau}^S_k\big|\mathcal{F}_k^s\right], ~\forall k\ge k_0.
\end{align}
On the other hand, by Proposition \ref{prop:vra}, 
\begin{align}\label{tilde-tau}
\mathbb{E}\left[\tilde{\tau}^S_k\right]&\le \left[\left(\frac{80L^2}{\tilde{\gamma} \eta_w}+30L^2\right)\gamma_k^2+\frac{2\gamma_k}{ \eta_w}+\frac{10\gamma_k^2}{\alpha_kn}\right]c_x\beta_{k} +\frac{8\lambda_k^2}{\tilde{\gamma} \eta_w}\left(n\sigma_\xi^2\right)+\frac{4\tilde{\gamma}n\sigma_\zeta^2\beta_k^2}{ \eta_w}+3n\sigma_\xi^2\lambda_k^2\notag\\
&\quad+\left[150L^2\alpha_k^2\tilde{\gamma}^2+\frac{4\tilde{\gamma}\left(\beta_k^2+20L^2\alpha_k^2\right)}{\eta_w}+\frac{10\alpha_k\tilde{\gamma}^2}{n}\right]c_s\beta_{k}\notag\\
&\le \left[\left(\frac{80L^2}{\tilde{\gamma} \eta_w}+30L^2\right)c_1^2\alpha_1+\frac{2c_1}{ \eta_w}+\frac{10c_1^2}{n}\right]c_xc_3\alpha_{k}^2 +\frac{8c_2\alpha_k^2}{\tilde{\gamma} \eta_w}\left(n\sigma_\xi^2\right)+\frac{4\tilde{\gamma}n\sigma_\zeta^2c_3^2\alpha_k^2}{ \eta_w}\notag\\
&\quad+3n\sigma_\xi^2c_2^2\alpha_k^2+\left[150L^2\alpha_1\tilde{\gamma}^2+\frac{4\tilde{\gamma}\left(c_3^2\alpha_1+20L^2\alpha_1\right)}{\eta_w}+\frac{10\tilde{\gamma}^2}{n}\right]c_sc_3\alpha_{k}^2\notag\\
&= C^S\alpha_{k}^2,
\end{align}
where 
\begin{align*}
C^S&:= \left[\left(\frac{80L^2}{\tilde{\gamma} \eta_w}+30L^2\right)c_1^2\alpha_1+\frac{2c_1}{ \eta_w}+\frac{10c_1^2}{n}\right]c_xc_3 +\frac{8c_2}{\tilde{\gamma} \eta_w}\left(n\sigma_\xi^2\right)+\frac{4\tilde{\gamma}n\sigma_\zeta^2c_3^2}{ \eta_w}\notag\\
&\quad+3n\sigma_\xi^2c_2^2+\left[150L^2\alpha_1\tilde{\gamma}^2+\frac{4\tilde{\gamma}\left(c_3^2\alpha_1+20L^2\alpha_1\right)}{\eta_w}+\frac{10\tilde{\gamma}^2}{n}\right]c_sc_3,
\end{align*}
the second inequality follows from (\ref{para-0}) and the nonincreasing property of $\alpha_{k}$. Taking expectation on both sides of (\ref{L-bound-4}), and substituting  (\ref{tilde-tau}) into it,
\begin{align}\label{L-bound-0}
	\mathbb{E}\left[\mathcal{L}^S_{k+1}\right]
	\le&\left(1-\frac{\delta^S}{2}\alpha_k\right)\mathbb{E}\left[\mathcal{L}^S_{k}\right]+C^S\alpha_k^2.
\end{align}
 Recursively, we have
 \begin{align}\label{L-sum}
 \mathbb{E}\left[\mathcal{L}^S_{k+1}\right]
 	\le&\Pi_{t=k_0}^k\left(1-\frac{\delta^S}{2}\alpha_t\right)\mathbb{E}\left[\mathcal{L}^S_1\right]+C^S\sum_{t= k_0}^k\Pi_{l=t+1}^k\left(1-\frac{\delta^S}{2}\alpha_l\right)\alpha_t^2, ~\forall k\ge k_0,
 \end{align}
where $\Pi_{l=t+1}^k\left(1-\frac{\delta^S}{2}\alpha_l\right)=1$ for $t=k$. When $\alpha_k\asymp\frac{1}{k^{b}}$ with $b\in(1/2,1)$, the first term on the right-hand side of (\ref{L-sum}) converges to zero exponentially. By \cite[Lemma I2]{Cardot2017fast}, the second term is of order $\mathcal{O}\left(\frac{1}{k^{b}}\right)$. Then,
\begin{equation*}
\mathbb{E}\left[\mathcal{L}^S_{k+1}\right]=\mathcal{O}\left(\frac{1}{k^{b}}\right).
\end{equation*}

\textbf{Part (ii).} When $\alpha_k=\frac{a_3}{k+a_1}$ with $a_1>0$, $
a_3> \frac{2}{\delta^S}
$, applying \cite[Lemma 4, Chapter 2]{polyak1987Introduction} to inequality (\ref{L-bound-0}) yields 
\begin{equation*}
	\mathbb{E}\left[\mathcal{L}^S_{k+1}\right]=\mathcal{O}\left(\frac{1}{k}\right).
\end{equation*} 
 The proof is complete.
\end{proof}
Theorem \ref{thm:x-opt-mse} establishes that VRA-DSGT attains an asymptotic convergence rate of $\mathcal{O}(k^{-1})$ in the mean square sense under diminishing stepsizes%\footnote{Equation (\ref{L-sum}) in the proof of Theorem \ref{thm:x-opt-mse} also provides a non-asymptotic upper bound for $\mathcal{L}^S_{k+1}$.}
, which matches the convergence rates of centralized stochastic gradient methods under perfect communication \cite{Nemi2009Robust,doan2020localSGD}.
	
 The following theorem establishes the convergence rate of VRA-DGT in the mean square sense under both constant and diminishing stepsizes.
\begin{thm}\label{thm:constant-mse}
	Suppose Assumptions \ref{ass:function}, \ref{ass:matrix} and \ref{ass:noise} (with $p_1=2$) hold. Let $\gamma_k = c_1 \alpha_k$,
	\begin{equation}\label{delta-1}
		\delta^D= \min\left\{\mu,\frac{\tilde{\gamma} \eta_w}{\alpha_1}-\frac{2}{c_1 \eta_w},c_1 \eta_w-\frac{8L^2}{n}\right\},~\tilde{\delta}^D=\left(\frac{10L^2}{ \tilde{\gamma}\eta_w}+30L^2\right)\max\left\{1,\left\|\mathbf{W}-\mathbf{I}\right\|^2c_1^2\right\},
	\end{equation}
{\small\begin{align}\label{tau-1}
		\tau^D=&\left[\left(\frac{10L^2}{ \tilde{\gamma}\eta_w}+30L^2\right)c_1^2+\frac{2c_1}{ \eta_w}+\frac{10c_1^2}{n}\right]c_x+10\tilde{\gamma}\left(\frac{\alpha_1L^2}{\eta_w}+\frac{\tilde{\gamma}}{n}\right)c_s,~\tilde{\tau}^D=\frac{4\tilde{\gamma}}{\eta_w}\left(c_s+n\sigma_\zeta^2\right),
\end{align}}where $c_1$ is a positive value, $c_s$ and $c_x$ are defined in Proposition \ref{prop:vra} and	$\eta_w$ is defined in (\ref{eta}). Then for VRA-DGT,
	\begin{itemize}
		\item[(i)]	$$\mathbb{E}\left[\mathcal{L}^D_{k+1}\right]=\mathcal{O}\left(\frac{1}{k^{b}}\right)$$
		if $\alpha_k\asymp\frac{1}{k^{b}}$ is nonincreasing,  $\beta_k \le c_2\alpha_k$ and that the parameters $b$,   $\tilde{\gamma}$, $c_1$ and $c_2$ satisfy $b\in(1/2,1)$, $c_2>0$ and
		\begin{equation}\label{cond-0}
			\tilde{\gamma} \in (0,1], ~\delta^D>0;
		\end{equation}
	\item[(ii)]	$$\mathbb{E}\left[\mathcal{L}^D_{k+1}\right]=\mathcal{O}\left(\frac{1}{k}\right)$$	if  $\alpha_k=\frac{a_2}{k+a_1}$, $\beta_k=\frac{c_2a_2}{k+a_1}$ and that the parameters $a_1$,  $a_2$, $c_1$, $c_2$ and $\tilde{\gamma}$ satisfy $a_1>0$, $a_2> \frac{2}{\delta^D}
	$, $c_2a_2\in[1,a_1+1)$ and the setting (\ref{cond-0});
		\item[(ii)]	\begin{equation}\label{non-asym}
			\mathbb{E}\left[\mathcal{L}^D_{k+1}\right]
			\le\left(1-\frac{\delta^D\tilde{\alpha}}{2}\right)^k\mathbb{E}\left[\mathcal{L}^D_1\right]+\frac{}{}\frac{4\tau^D}{\delta^D}\beta_{k}+\frac{4\tilde{\tau}^D}{\delta^D\tilde{\alpha}}\beta_{k}^2,
		\end{equation}	if $\alpha_k\equiv \tilde{\alpha}$, $\beta_k\asymp \frac{1}{k^b}$ or $\beta_k=\frac{a_2}{k+a_1}$,  and that the  parameters $b$, $a_1$, $a_2$, $c_1$, $\tilde{\gamma}$ and  $ \tilde{\alpha}$ satisfy
		{\footnotesize\begin{align}
			&b\in \left(\frac{1}{2},1\right), a_1\ge 0, a_2\in[1,a_1+1),\frac{4-3\delta^D\tilde{\alpha}}{4-2\delta^D\tilde{\alpha}}\ge \sup_{k\ge 1}\frac{\beta_{k}^2}{\beta_{k+1}^2}, \tilde{\alpha} \le  \min \left\{1,\frac{1}{\left(\frac{10L^2}{ \tilde{\gamma}\eta_w}+30L^2\right)n},\frac{2\tilde{\delta}^D}{\delta^D}\right\},\label{para-2}
			\end{align}}
		 and the setting (\ref{cond-0}).
		\end{itemize} 
\end{thm}
	\begin{proof}
	\textbf{Parts (i) and (ii)}. 
	Note that the difference between  $\mathcal{L}^D_{k+1}$ and $\mathcal{L}^S_{k+1}$ is the  term  $\left\|\mathbf{e}_k^q\right\|^2$.
	According to the similar analysis of  (\ref{L-bound-2}), 
	\begin{align}\label{L-bound-GT-0}
		&\mathbb{E}\left[\mathcal{L}^D_{k+1}\right]\notag\\
		\le& (1-\mu\alpha_k)(f(\bar{x}_k)-f^*)+\left[1-\tilde{\gamma} \eta_w+\frac{2\alpha_k^2}{\gamma_k \eta_w}+\left(\frac{10L^2}{\tilde{\gamma} \eta_w}+30L^2\right)\alpha_k^2\right]\mathbb{E}\left[\left\|\mathbf{y}_{k}-\bar{\mathbf{y}}_{k}\right\|^2\right]\notag\\
		&+\left[1-\gamma_k \eta_w+\frac{8\alpha_kL^2}{n}+\left(\frac{10L^2}{\tilde{\gamma} \eta_w}+30L^2\right)\left\|\mathbf{W}-\mathbf{I}\right\|^2\gamma_k^2\right]\mathbb{E}\left[\left\|\mathbf{x}_k-\bar{\mathbf{x}}_k\right\|^2\right]\notag\\
		&+\left[\left(\frac{10L^2}{\tilde{\gamma} \eta_w}+30L^2\right)\gamma_k^2+\frac{2\gamma_k}{ \eta_w}+\frac{10\gamma_k^2}{\alpha_kn}\right]\mathbb{E}\left[\left\|\mathbf{e}_k^x\right\|^2\right]\notag\\
		&+\left[\frac{4\tilde{\gamma}\beta_k^2+10\tilde{\gamma}L^2\alpha_k^2}{\eta_w}+\frac{10\alpha_k\tilde{\gamma}^2}{n}\right]\mathbb{E}\left[\left\|\mathbf{e}_k^s\right\|^2\right]\notag\\
		& -\left[\alpha_k-\left(\frac{10L^2}{\tilde{\gamma} \eta_w}+30L^2\right)n\alpha_k^2\right]\mathbb{E}\left[\|\bar{q}_{k}\|_2^2\right] +\frac{4\tilde{\gamma}\beta_k^2}{ \eta_w}n\sigma_\zeta^2\notag\\
		&\le \left(1-\delta^D\alpha_k+\tilde{\delta}^D\alpha_k^2\right)\mathbb{E}\left[\mathcal{L}^D_k\right]+\tau^D\alpha_{k}\beta_{k}+\tilde{\tau}^D\beta_k^2,
	\end{align}
	where  $\delta^D$, $\tilde{\delta}^D$, $\tau^D$ and $\tilde{\tau}^D$  are defined in  (\ref{delta-1}) and (\ref{tau-1}), the second inequality follows from Proposition \ref{prop:vra}.
Then there exists a large enough integer $k_0$ such that
\begin{align}\label{L-sum-0}
	\mathbb{E}\left[\mathcal{L}^D_{k+1}\right]
	\le&\left(1-\frac{\delta^D}{2}\alpha_k\right)\mathbb{E}\left[\mathcal{L}^D_{k}\right]+\left(\tau^Dc_2+\tilde{\tau}^Dc_2^2\right)\alpha_k^2,~\forall k\ge k_0.
\end{align}
 Recursively, we have
{\small\begin{align}\label{L-sum-1}
	\mathbb{E}\left[\mathcal{L}^D_{k+1}\right]
	\le&\Pi_{t=k_0}^k\left(1-\frac{\delta^D}{2}\alpha_t\right)\mathbb{E}\left[\mathcal{L}^D_1\right]+\left(\tau^Dc_2+\tilde{\tau}^Dc_2^2\right)\sum_{t= k_0}^k\Pi_{l=t+1}^k\left(1-\frac{\delta^D}{2}\alpha_l\right)\alpha_t^2, ~\forall k\ge k_0.
\end{align}}Combining  (\ref{L-sum-0}) with \cite[Lemma I2]{Cardot2017fast} and  (\ref{L-sum-1}) with \cite[Lemma 4, Chapter 2]{polyak1987Introduction} yields the desired results of Parts (i) and (ii) respectively.

\textbf{Part (iii)}.	Given (\ref{L-bound-GT-0}) and the fact $\alpha_k \equiv \tilde{\alpha}$, $\tilde{\alpha} \le \frac{2\tilde{\delta}^D}{\delta^D}$, we get
\begin{align*}
		\mathbb{E}\left[\mathcal{L}^D_{k+1}\right]
		\le&\left(1-\frac{\delta^D\tilde{\alpha}}{2}\right)\mathbb{E}\left[\mathcal{L}^D_k\right]+\tau^D\tilde{\alpha}\beta_{k}+\tilde{\tau}^D\beta_k^2\\
		\le&\left(1-\frac{\delta^D\tilde{\alpha}}{2}\right)^k\mathbb{E}\left[\mathcal{L}^D_1\right]+\tau^D\tilde{\alpha}\sum_{t=1}^k\left(1-\frac{\delta^D\tilde{\alpha}}{2}\right)^{k-t}\beta_{t}+\tilde{\tau}^D\sum_{t=1}^k\left(1-\frac{\delta^D\tilde{\alpha}}{2}\right)^{k-t}\beta_{t}^2.	
	\end{align*}
	 By the fact $\sup_{k\ge 1}\frac{\beta_{k}^2}{\beta_{k+1}^2}\le 1+\frac{\delta^D\tilde{\alpha}}{4-2\delta^D\tilde{\alpha}}$ and  \cite[Lemma 7]{Li2022Perform}, the last two terms on the right hand side of above inequality satisfy $$\tau^D\tilde{\alpha}\sum_{t=1}^k\left(1-\frac{\delta^D\tilde{\alpha}}{2}\right)^{k-t}\beta_{t}\le \frac{4\tau^D}{\delta^D}\beta_{k},~\tilde{\tau}^D \sum_{t=1}^k\left(1-\frac{\delta^D\tilde{\alpha}}{2}\right)^{k-t}\beta_{t}^2\le \frac{4\tilde{\tau}^D}{\delta^D\tilde{\alpha}}\beta_{k}^2.$$ Then, we arrive at (\ref{non-asym}).
	The proof is complete.
	\end{proof}

		Theorem \ref{thm:constant-mse} provides both asymptotic and non-asymptotic convergence rates for VRA-DGT  in the mean square sense under diminishing and constant stepsizes. By choosing $\beta_k = \mathcal{O}(k^{-1})$, $\alpha_k=\mathcal{O}(k^{-1})$ or $\alpha_k\equiv \tilde{\alpha}$, VRA-DGT achieves the convergence rate  $\mathcal{O}(k^{-1})$, which  is faster than  $\mathcal{O}(k^{-1/2})$ of DGD \cite{Reisizadeh2023Vary} and  $\mathcal{O}(k^{-1+\delta})$  (for any $\delta > 0$) of VRA-GT \cite{zhao2023VRA}. This convergence rate also complements the results in \cite{pu2020robust} that establish convergence to a neighborhood of the optimum under the constant stepsize.

\subsection{Convergence Rate in the Almost Sure Sense}

In this subsection, we provide the convergence rate of VRA-D(S)GT in the almost sure sense.  Different from the convergence
	guarantees in expectation, the almost sure convergence rate describes the behavior of individual sample paths generated by stochastic algorithms,  which  corresponds to the actual realizations of the stochastic algorithms used in practice \cite{Liu2022asr}. For non-distributed stochastic optimization problems, 
	the almost sure convergence rate of stochastic gradient algorithms has been widely studied  \cite{Peggy2025Newton, Liu2022asr, chen2006stochastic}. 
	For distributed stochastic optimization problems, Xin et al. \cite{Xin2012improv} show that the distributed stochastic gradient tracking method achieves the  convergence rate of $o(k^{-1+\delta})$ in the almost sure sense, where  $\delta > 0$ is arbitrarily small. More recently, Deng et al. \cite{Deng2025Almost} show that a stochastic
	unified decentralized algorithm achieves  the convergence rate of $\mathcal{O}\left(\frac{\ln(k)}{k^b}\right)$, $\forall b\in(0.5,1)$ in the almost sure sense. However, for imperfect-communication networks, existing studies \cite{Doan2021Quanti,Sri2011async,Zhang2019Sign,zhao2023VRA} have primarily focused on almost sure convergence. In what follows, we complement these results by establishing the almost sure convergence rate of VRA-D(S)GT.
\begin{thm}\label{thm:as-DSGT}
		Suppose  Assumptions \ref{ass:function}, \ref{ass:matrix}, \ref{ass:noise} (with $p_1>2$) and \ref{ass:stochastic gradient} (with $p_2>2$) hold and that $L_i(\xi_i)<\hat{c}$, $\forall i\in\mathcal{V}$ almost surely, where $L_i(\xi_i)$ is defined in Assumption \ref{ass:stochastic gradient} and $\hat{c}$ is some positive constant. Suppose also that $\alpha_k\asymp\frac{1}{k^{b}}$ is nonincreasing, 
		$$\gamma_k=c_1\alpha_k,\lambda_k=c_2\alpha_k,\beta_k\le c_3 \alpha_k,
		$$
		where parameters  $\tilde{\gamma}$, $b$, $c_1$, $c_2$  and $c_3$ satisfy 
		\begin{align*}
			b\in \left(\frac{p_2}{2(p_2-2)},1\right),~1+\frac{\tilde{\gamma} \eta_w}{4-2\tilde{\gamma} \eta_w}\ge \sup_{k\ge 1}\frac{\alpha_k^2}{\alpha_{k+1}^2},~c_1<\frac{1}{\alpha_1},~c_2>\max\{10(L^2+\hat{c}^2),16\}
		\end{align*}
		and the setting (\ref{para}), $\eta_w$ is defined in (\ref{eta}).
		Then for VRA-DSGT,
		$$\mathcal{L}^S_{k+1}=\mathcal{O}\left(\alpha_k\ln(k+1)\right)~~~\text{a.s.}$$
	\end{thm}
	\begin{proof} By the definition of $\mathcal{L}^S_{k+1}$ in (\ref{L-fun-2}), we finish the proof by deriving individual bounds for the gradient tracking error $\left\|\mathbf{y}_{k+1}-\bar{\mathbf{y}}_{k+1}\right\|^2$,  the agreement error $\left\|\mathbf{x}_{k+1}-\bar{\mathbf{x}}_{k+1}\right\|^2$,  the optimality gap  $f(\bar{x}_{k+1}) - f^*$ with  VRA estimation error $\left\|\mathbf{e}_{k+1}^q\right\|^2$ respectively.

	\textbf{$\left\|\mathbf{y}_{k+1}-\bar{\mathbf{y}}_{k+1}\right\|^2$.} By the definition of  $\bar{\mathbf{y}}_{k+1}$ and the recursions (\ref{alg:q}) and (\ref{e-s}),
		\begin{align*}
			\mathbf{y}_{k+1}-\bar{\mathbf{y}}_{k+1}
			&=\left(\mathbf{W}_{\tilde{\gamma}}-\frac{\mathbf{1}\mathbf{1}^\intercal}{n}\right)(\mathbf{y}_{k}-\bar{\mathbf{y}}_{k})+\left(\mathbf{I}-\frac{\mathbf{1}\mathbf{1}^\intercal}{n}\right)\left(\mathbf{q}_{k+1}-\mathbf{q}_{k}+\tilde{\gamma}\mathbf{e}_{k+1}^s-\tilde{\gamma}\mathbf{e}_{k}^s\right)\\
			&=\left(\mathbf{W}_{\tilde{\gamma}}-\frac{\mathbf{1}\mathbf{1}^\intercal}{n}\right)(\mathbf{y}_{k}-\bar{\mathbf{y}}_{k})+\left(\mathbf{I}-\frac{\mathbf{1}\mathbf{1}^\intercal}{n}\right)\left(-\lambda_k\underbrace{\left(\mathbf{q}_{k}-\mathbf{g}_{k}\right)}_{\mathbf{e}_k^q}+\lambda_k\underbrace{\left(\tilde{\mathbf{g}}_{k+1}^{'}-\mathbf{g}_{k+1}\right)}_{\epsilon_k}\right.\\
			&\quad\left.+\underbrace{\lambda_k\left(\mathbf{g}_{k+1}-\mathbf{g}_{k}\right)+(1-\lambda_k)\left(\tilde{\mathbf{g}}_{k+1}^{'}-\tilde{\mathbf{g}}_k\right)}_{\tilde{\epsilon}_k}-\tilde{\gamma}\beta_k\mathbf{e}_{k}^s+\tilde{\gamma}\beta_k\mathbf{W}\zeta_{k+1}^s\right).
		\end{align*}
		Then,
		\begin{align*}%\label{y-bary}
			\|\mathbf{y}_{k+1}-\bar{\mathbf{y}}_{k+1}\|
			&=\left\|\left(\mathbf{W}_{\tilde{\gamma}}-\frac{\mathbf{1}\mathbf{1}^\intercal}{n}\right)^k(\mathbf{y}_1-\bar{\mathbf{y}}_1)+\sum_{t= 1}^k\left(\mathbf{W}_{\tilde{\gamma}}-\frac{\mathbf{1}\mathbf{1}^\intercal}{n}\right)^{k-t}\left(-\lambda_t\mathbf{e}^q_t-\tilde{\gamma}\beta_t\mathbf{e}_{t}^s\right)\right.\notag\\
			&\left.\quad+\sum_{t= 1}^k\left(\mathbf{W}_{\tilde{\gamma}}-\frac{\mathbf{1}\mathbf{1}^\intercal}{n}\right)^{k-t}\left(\lambda_t\epsilon_t+\tilde{\gamma}\beta_t\mathbf{W}\zeta_{t+1}^s\right)+\sum_{t= 1}^k\left(\mathbf{W}_{\tilde{\gamma}}-\frac{\mathbf{1}\mathbf{1}^\intercal}{n}\right)^{k-t}\tilde{\epsilon}_t\right\|\notag\\
			&\le \left(1-\tilde{\gamma}\eta_w\right)^k\|\mathbf{y}_1-\bar{\mathbf{y}}_1\|+\underbrace{\sum_{t= 1}^k\left(1-\tilde{\gamma}\eta_w\right)^{k-t}\left(\lambda_t\|\mathbf{e}^q_t\|+\tilde{\gamma}\beta_t\|\mathbf{e}_{t}^s\|\right)}_{P_1}\notag\\
			&\quad+\underbrace{\left\|\sum_{t= 1}^k\left(\mathbf{W}_{\tilde{\gamma}}-\frac{\mathbf{1}\mathbf{1}^\intercal}{n}\right)^{k-t}\left(\lambda_t\epsilon_t+\tilde{\gamma}\beta_t\mathbf{W}\zeta_{t+1}^s\right)\right\|}_{P_2}+\underbrace{\sum_{t= 1}^k\left(1-\tilde{\gamma}\eta_w\right)^{k-t}\|\tilde{\epsilon}_t\|}_{P_3}.
		\end{align*}
		On the right hand side of the above inequality, the first term converges exponentially fast to 0. We next bound the terms $P_1$, $P_2$ and $P_3$ respectively. 
		
		Lemma \ref{lem:e-o} in Appendix B and Proposition \ref{prop:vra-1} imply that $\|\mathbf{e}^q_t\|^2$ is finite and
		$\left\|\mathbf{e}^s_{t}\right\|^2\le \tilde{c}_s\beta_t\ln (t+1)$ almost surely.
		Then, there exists a positive finite random variable $C_1$ such that
		\begin{align*}
			P_1\le \sum_{t= 1}^k\left(1-\tilde{\gamma}\eta_w\right)^{k-t} C_1\left(\frac{1}{t^b}+\frac{1}{t^{3b/2}}\sqrt{\ln (t+1)}\right)=\mathcal{O}\left(\frac{1}{t^b}\right)~~~\text{a.s.}.
		\end{align*}
		
		Note that $\epsilon_t$ and $\zeta_{t}^s$ are martingale sequences and $\lambda_k=c_2\alpha_k,\beta_k\asymp \alpha_k$. Then by a similar analysis of Proposition \ref{prop:vra-1}, we have {\small$P_2=\mathcal{O}\left(\sqrt{\alpha_k\ln (k+1)}\right)$} almost surely.

		By the Lipschitz continuous of $\nabla f_j(x)$ and $\nabla F_j(x;\xi_j),j\in\mathcal{V}$, we have
		{\small\begin{align*}
			&\|\tilde{\epsilon}_t\|\\
			&=\left\|\lambda_k\left(\mathbf{g}_{k+1}-\mathbf{g}_{k}\right)+(1-\lambda_k)\left(\tilde{\mathbf{g}}_{k+1}^{'}-\tilde{\mathbf{g}}_k\right)\right\|\notag\\
			&\le \left(\lambda_kL+(1-\lambda_k)\max_{i\in\mathcal{V}}L_i(\xi_i)\right)\|\mathbf{x}_{k+1}-\mathbf{x}_k\|\notag\\
			&\le \left(\lambda_kL+(1-\lambda_k)\hat{c}\right)\left(\gamma_k\|\mathbf{W}-\mathbf{I}\|\left\|\mathbf{x}_k-\bar{\mathbf{x}}_k\right\|
			+\gamma_k\left\|\mathbf{e}_k^x\right\|
			+\alpha_k\left\|\mathbf{y}_k-\bar{\mathbf{y}}_k\right\|
			+\alpha_k\tilde{\gamma}\left\|\mathbf{e}^s_{k}\right\|+\alpha_k\sqrt{n}\left\|\bar{q}_{k}\right\|_2\right),
		\end{align*}
		}where $\bar{q}_{k}:=\frac{1}{n} \sum_{j=1}^nq_{j,k}$, {\small$\bar{e}_{k}^s:=\frac{1}{n} \sum_{j=1}^n\left(z_{j,k}^s-\sum_{j\in \mathcal{N}_{i}}s_{j,k}\right)$}, the second inequality follows from the condition $L_i(\xi_i)<\hat{c}$, $\forall i\in\mathcal{V}$ almost surely and  (\ref{x-bound}) in Appendix A.
		For the last term  on the right hand side of the above inequality,
		\begin{align*}
			\alpha_k\sqrt{n}\left\|\bar{q}_k\right\|_2^2=&\alpha_k\sqrt{n}\left\|\bar{e}_k^q+\frac{1}{n}\sum_{j=1}^n\nabla f_j(x_{j,k})-\nabla f(\bar{x}_k)+\nabla f(\bar{x}_k)\right\|_2\notag\\
			\le&\alpha_k\sqrt{n}\left\|\bar{e}_k^q\right\|_2+\alpha_k\sqrt{n}\left\|\frac{1}{n}\sum_{j=1}^n\nabla f_j(x_{j,k})-\nabla f(\bar{x}_k)\right\|_2+\alpha_k\sqrt{n}\left\|\nabla f(\bar{x}_k)\right\|_2\notag\\
			\le&\alpha_k \left\|\mathbf{e}_k^q\right\|+\alpha_kL \left\|\mathbf{x}_k-\mathbf{\bar{x}}_k\right\|+\alpha_k\sqrt{n}\left\|\nabla f(\bar{x}_k)\right\|_2,
		\end{align*}
		where $\bar{e}_k^q:=\bar{q}_{k}-\frac{1}{n}\sum_{j=1}^n\nabla f(x_{j,k})$.
		Then,
		\begin{align}\label{epsilon}
			\|\tilde{\epsilon}_k\|
			&\le \left(\lambda_kL+(1-\lambda_k)\hat{c}\right) \left(\left(\gamma_k\|\mathbf{W}-\mathbf{I}\|+\alpha_kL\right)\left\|\mathbf{x}_k-\bar{\mathbf{x}}_k\right\|
			+\gamma_k\left\|\mathbf{e}_k^x\right\|
			+\alpha_k\left\|\mathbf{y}_k-\bar{\mathbf{y}}_k\right\|\right.
			\notag\\
			&\quad\left.+\alpha_k\left\|\mathbf{e}_k^q\right\|
			+\alpha_k\tilde{\gamma}\left\|\mathbf{e}^s_{k}\right\|+\alpha_k\sqrt{n}\|\nabla f(\bar{x}_k)\|_2\right)\notag\\
			&\le \left(\lambda_kL+(1-\lambda_k)\hat{c}\right) \left(\left(\gamma_k\|\mathbf{W}-\mathbf{I}\|+\alpha_kL\right)\left\|\mathbf{x}_k-\bar{\mathbf{x}}_k\right\|
			+\gamma_k\tilde{c}_x\beta_k\ln (k+1)
			\right.
			\notag\\
			&\quad\left.+\alpha_k\left\|\mathbf{y}_k-\bar{\mathbf{y}}_k\right\|+\alpha_k\left\|\mathbf{e}_k^q\right\|
			+\alpha_k\tilde{\gamma}\tilde{c}_s\beta_k\ln (k+1)+\alpha_k\sqrt{n}\sqrt{2L(f(\bar{x}_k)-f^*)}\right),
		\end{align}
		where the second inequality follows from Proposition \ref{prop:vra-1} and the Lipschitz smoothness of $f(x)$. On the other hand, Lemma \ref{lem:e-o} in Appendix B implies the terms $\left\|\mathbf{x}_k-\bar{\mathbf{x}}_k\right\|$, $\left\|\mathbf{y}_k-\bar{\mathbf{y}}_k\right\|$, $\left\|\mathbf{e}_k^q\right\|$ and $f(\bar{x}_k)-f^*$ are finite almost surely. Then, $\|\tilde{\epsilon}_t\|=\mathcal{O}\left(\frac{1}{t^b}\right)$, which implies $P_3=\mathcal{O}\left(\frac{1}{t^b}\right)$ almost surely. 
		
		Combining the bounds of $P_1$, $P_2$ and $P_3$, we get that
		\begin{equation}\label{y-bound}
			\left\|\mathbf{y}_{k+1}-\bar{\mathbf{y}}_{k+1}\right\|^2=\mathcal{O}\left(\alpha_k\ln(k+1)\right).
		\end{equation}
		
		\textbf{$\left\|\mathbf{x}_{k+1}-\bar{\mathbf{x}}_{k+1}\right\|^2$.}  By the definition of  $\bar{\mathbf{x}}_{k+1}$ and the recursion (\ref{alg:x}),
		\begin{align*}
			\|\mathbf{x}_{k+1}-\bar{\mathbf{x}}_{k+1}\|
			&=\left\|\left(\mathbf{W}_{\gamma_k}-\frac{\mathbf{1}\mathbf{1}^\intercal}{n}\right)(\mathbf{x}_{k}-\bar{\mathbf{x}}_{k})-\alpha_k\left(\mathbf{y}_{k}-\bar{\mathbf{y}}_{k}\right)+\gamma_k\left(\mathbf{I}-\frac{\mathbf{1}\mathbf{1}^\intercal}{n}\right)\mathbf{e}_{k}^x\right\|\notag\\
			&\le \left(1-\gamma_k\eta_w\right)\|\mathbf{x}_{k}-\bar{\mathbf{x}}_{k}\|+\alpha_k\|\mathbf{y}_{k}-\bar{\mathbf{y}}_{k}\|+\gamma_k\|\mathbf{e}_{k}^x\|.
		\end{align*}
		Recursively, we have
		\begin{align}\label{x-bound-1}
			\|\mathbf{x}_{k+1}-\bar{\mathbf{x}}_{k+1}\|
			&\le \Pi_{t=1}^k\left(1-\gamma_t\eta_w\right)\|\mathbf{x}_{1}-\bar{\mathbf{x}}_{1}\|+\sum_{t= 1}^k\left(\Pi_{l=t+1}^k\left(1-\gamma_t\eta_w\right)\right)\alpha_t\|\mathbf{y}_{t}-\bar{\mathbf{y}}_{t}\|\notag\\
			&\quad+\sum_{t= 1}^k\left(\Pi_{l=t+1}^k\left(1-\gamma_t\eta_w\right)\right)\gamma_t\|\mathbf{e}_{t}^x\|.
		\end{align}
		The first term on the right hand side of (\ref{x-bound-1}) converges exponentially fast to 0.  By (\ref{y-bound}) and \cite[Lemma I2]{Cardot2017fast}, the second term on the right hand side of (\ref{x-bound-1})
		\begin{equation*}
			\sum_{t= 1}^k\left(\Pi_{l=t+1}^k\left(1-\gamma_t\eta_w\right)\right)\alpha_t\|\mathbf{y}_{t}-\bar{\mathbf{y}}_{t}\|=\mathcal{O}\left(\sqrt{\alpha_k\ln(k+1)}\right)~~~\text{a.s.}.
		\end{equation*}
		For the last term on the right hand side of (\ref{x-bound-1}),
		{\small$$\sum_{t= 1}^k\left(\Pi_{l=t+1}^k\left(1-\gamma_t\eta_w\right)\right)\gamma_t\|\mathbf{e}_{t}^x\|\le \sum_{t= 1}^k\left(\Pi_{l=t+1}^k\left(1-\gamma_t\eta_w\right)\right)\gamma_t\tilde{c}_x\beta_t\ln (t+1)=\mathcal{O}\left(\alpha_k\ln(k+1)\right)~~~\text{a.s.},$$}where the inequality follows from Proposition \ref{prop:vra-1}, and the equality holds by the condition $\beta_k\asymp \alpha_{k}$ and \cite[Lemma I2]{Cardot2017fast}.
		Hence, $\left\|\mathbf{x}_{k+1}-\bar{\mathbf{x}}_{k+1}\right\|^2=\mathcal{O}\left(\alpha_k\ln(k+1)\right)$ almost surely.

		\textbf{$f(\bar{x}_{k+1})-f^*+\left\|\mathbf{e}_{k+1}^q\right\|_2^2$.} Recall (\ref{f-bound}) in Appendix A.4, 
		\begin{align*}
			f(\bar{x}_{k+1})
			&\le f(\bar{x}_k)-\frac{\alpha_k}{2}\|\nabla f(\bar{x}_k)\|_2^2-\alpha_k\|\bar{q}_{k}\|_2^2+4\alpha_k\|\bar{q}_{k}-\nabla f(\bar{x}_k)\|_2^2+5\alpha_k\left\|\tilde{\gamma}\bar{e}^s_k+\frac{\gamma_k}{\alpha_k}\bar{e}^x_{k}\right\|_2^2,
		\end{align*}
		where {\small$\bar{e}_{k}^s:=\frac{1}{n} \sum_{j=1}^n\left(z_{j,k}^s-\sum_{j\in \mathcal{N}_{i}}s_{j,k}\right)$ and $\bar{e}_{k}^x:=\frac{1}{n} \sum_{j=1}^n\left(z_{j,k}^x-\sum_{j\in \mathcal{N}_{i}}x_{j,k}\right)$.}
		For the last two terms on the right hand side of the above inequality, we have
		\begin{align*}
			4\alpha_k\|\bar{q}_{k}-\nabla f(\bar{x}_k)\|_2^2=&4\alpha_k\left\|\bar{q}_{k}-\frac{1}{n}\sum_{j=1}^n\nabla f(x_{j,k})+\frac{1}{n}\sum_{j=1}^n\nabla f_j(x_{j,k})-\nabla f(\bar{x}_k)\right\|_2^2\\
			%\le&8\alpha_k\left\|\bar{e}_{k}^q\right\|^2+8\alpha_k\left\|\frac{1}{n}\sum_{j=1}^n\nabla f_j(x_{j,k})-\nabla f(\bar{x}_k)\right\|_2^2\\
			\le&\frac{8\alpha_k}{n}\left\|\mathbf{e}_{k}^q\right\|^2+\frac{8\alpha_kL^2}{n}\left\|\mathbf{x}_k-\bar{\mathbf{x}}_k\right\|^2,
		\end{align*}
		and
		\begin{align*}
			&5\alpha_k\left\|\tilde{\gamma}\bar{e}^s_k+\frac{\gamma_k}{\alpha_k}\bar{e}^x_{k}\right\|_2^2\le \frac{10\alpha_k\tilde{\gamma}^2}{n}\|\mathbf{e}^s_{k}\|^2+\frac{10\gamma_k^2}{\alpha_kn}\|\mathbf{e}^x_{k}\|^2.
		\end{align*}
		Then,
		\begin{align*}
			f(\bar{x}_{k+1})
			&\le f(\bar{x}_k)-\frac{\alpha_k}{2}\|\nabla f(\bar{x}_k)\|_2^2-\alpha_k\|\bar{q}_{k}\|_2^2+\frac{8\alpha_k}{n}\left\|\mathbf{e}_{k}^q\right\|^2+\frac{8\alpha_kL^2}{n}\left\|\mathbf{x}_k-\bar{\mathbf{x}}_k\right\|^2\notag\\
			&\quad+\frac{10\alpha_k\tilde{\gamma}^2}{n}\|\mathbf{e}^s_{k}\|^2+\frac{10\gamma_k^2}{\alpha_kn}\|\mathbf{e}^x_{k}\|^2\notag\\
			&\le f(\bar{x}_k)-\frac{\alpha_k}{2}\|\nabla f(\bar{x}_k)\|_2^2-\alpha_k\|\bar{q}_{k}\|_2^2+\frac{16\alpha_k}{n}\left\|\mathbf{e}_{k}^q+m_k\right\|^2+\frac{16\alpha_k}{n}\left\|m_k\right\|^2\notag\\
			&\quad+\frac{8\alpha_kL^2}{n}\left\|\mathbf{x}_k-\bar{\mathbf{x}}_k\right\|^2+\frac{10\alpha_k\tilde{\gamma}^2}{n}\|\mathbf{e}^s_{k}\|^2+\frac{10\gamma_k^2}{\alpha_kn}\|\mathbf{e}^x_{k}\|^2,
		\end{align*}
		where 
		\begin{align*}
			m_k:=\sum_{t= 1}^{k-1}\left(\Pi_{l=t+1}^{k-1}(1-\lambda_l)\right)\lambda_t\left(\tilde{\mathbf{g}}_{t+1}^{'}-\mathbf{g}_{t+1}\right)~ \forall k\ge 2,~m_1=\mathbf{0}.
		\end{align*}
		Subtract $f^*$ from both sides of the above inequality, we have by the strong convexity of $f(x)$ that
		\begin{align}\label{f-bound-1}
			f(\bar{x}_{k+1})-f^*
			&\le \left(1-\mu\alpha_k\right)\left(f(\bar{x}_k)-f^*\right)-\alpha_k\|\bar{q}_{k}\|_2^2+\frac{16\alpha_k}{n}\left\|\mathbf{e}_k^q+m_k\right\|^2+\frac{16\alpha_k}{n}\left\|m_k\right\|^2\notag\\
			&\quad+\frac{8\alpha_kL^2}{n}\left\|\mathbf{x}_k-\bar{\mathbf{x}}_k\right\|^2+\frac{10\alpha_k\tilde{\gamma}^2}{n}\|\mathbf{e}^s_{k}\|^2+\frac{10\gamma_k^2}{\alpha_kn}\|\mathbf{e}^x_{k}\|^2.
		\end{align}
		By the definition of $\mathbf{e}^q_k$ and the recursion (\ref{alg:q}),
		\begin{align*}
			\mathbf{e}_{k+1}^q %&=(1-\lambda_k)\mathbf{e}_k^q+(1-\lambda_k)\left(\mathbf{g}_k-\tilde{\mathbf{g}}_k\right)+\left(\tilde{\mathbf{g}}_{k+1}^{'}-\mathbf{g}_{k+1}\right)\\
			&=(1-\lambda_k)\mathbf{e}_k^q+(1-\lambda_k)\underbrace{\left(\mathbf{g}_k-\mathbf{g}_{k+1}+\tilde{\mathbf{g}}_{k+1}^{'}-\tilde{\mathbf{g}}_k\right)}_{\hat{\epsilon}_k}+\lambda_k\underbrace{\left(\tilde{\mathbf{g}}_{k+1}^{'}-\mathbf{g}_{k+1}\right)}_{\epsilon_k}
			%\\
			%&=\Pi_{t=1}^k(1-\lambda_t)\mathbf{e}_1^q+\sum_{t= 1}^k\left(\Pi_{l=t+1}^k(1-\lambda_t)\right)\hat{\epsilon}_t+\sum_{t= 1}^k\left(\Pi_{l=t+1}^k(1-\lambda_t)\right)\lambda_t\epsilon_t
			.
		\end{align*}
		Then,
		\begin{align}\label{bar-e-bound}
			&\left\|\mathbf{e}_{k+1}^q+m_{k+1}\right\|^2\notag\\
			%&=\left\|\mathbf{e}_{k+1}^q-\sum_{t= 1}^k\left(\Pi_{l=t+1}^k(1-\lambda_t)\right)\lambda_t\epsilon_t\right\|^2\notag\\
			&=\left\|(1-\lambda_k)\mathbf{e}_k^q+(1-\lambda_k)\hat{\epsilon}_k+\lambda_k\epsilon_k-\sum_{t= 1}^k\left(\Pi_{l=t+1}^k(1-\lambda_t)\right)\lambda_t\epsilon_t\right\|^2\notag\\
			&=\left\|(1-\lambda_k)\left(\mathbf{e}_k^q+m_k\right)+(1-\lambda_k)\hat{\epsilon}_k\right\|^2\notag\\
			&\le  (1-\lambda_k)\left\|\mathbf{e}_k^q+m_k\right\|^2+\frac{1}{\lambda_k}\left\|\hat{\epsilon}_k\right\|^2\notag\\
			&\le  (1-\lambda_k)\left\|\mathbf{e}_k^q+m_k\right\|^2+\frac{2\left(L^2+\max_{i\in\mathcal{V}}L_i^2(\xi_i)\right)}{\lambda_k n}\left\|\mathbf{x}_{k+1}-\mathbf{x}_k\right\|^2\notag\\
			&\le  (1-\lambda_k)\left\|\mathbf{e}_k^q+m_k\right\|^2+\frac{2(L^2+\max_{i\in\mathcal{V}}L_i^2(\xi_i))}{\lambda_k n}\left(5\alpha_k^2n\left\|\bar{q}_k\right\|_2^2+5\gamma_k^2\left\|\mathbf{W}-\mathbf{I}\right\|^2\left\|\mathbf{x}_k-\bar{\mathbf{x}}_k\right\|^2\right.\notag\\
			&\quad\left.+5\gamma_k^2\left\|\mathbf{e}_k^x\right\|^2+5\alpha_k^2\left\|\mathbf{y}_k-\bar{\mathbf{y}}_k\right\|^2+5\alpha_k^2\left\|\mathbf{e}_k^s\right\|^2\right),
		\end{align}
		where 
		the second inequality  follows from the Lipschitz continuous of $\nabla f_j(x)$ and $\nabla F_j(x;\xi_j),j\in\mathcal{V}$, the last inequality follows from (\ref{x-bound}). Combining (\ref{bar-e-bound}) with (\ref{f-bound-1}) implies
	{\small	\begin{align*}
			&f(\bar{x}_{k+1})-f^*+\left\|\mathbf{e}_{k+1}^q+m_{k+1}\right\|^2\\
			&\le \left(1-\mu\alpha_k\right)\left(f(\bar{x}_k)-f^*\right)+\left(1-\lambda_k+16\alpha_k\right)\left\|\mathbf{e}_k^q+m_k\right\|^2-\alpha_k\left(1-\frac{10(L^2+\max_{i\in\mathcal{V}}L_i^2(\xi_i))}{c_2}\right)\|\bar{q}_{k}\|_2^2\notag\\
			&\quad+16\alpha_k\left\|m_k\right\|^2+\alpha_k\left(\frac{8L^2}{n}+\frac{10(L^2+\max_{i\in\mathcal{V}}L_i^2(\xi_i))c_2^2}{c_1 n}\left\|\mathbf{W}-\mathbf{I}\right\|^2\right)\left\|\mathbf{x}_k-\bar{\mathbf{x}}_k\right\|^2\\
			&\quad+\alpha_k\left(\frac{10\tilde{\gamma}^2}{n}+\frac{10(L^2+\max_{i\in\mathcal{V}}L_i^2(\xi_i))}{c_1 n}\right)\|\mathbf{e}^s_{k}\|^2+\alpha_k\left(\frac{10c_1^2}{n}+\frac{10(L^2+\max_{i\in\mathcal{V}}L_i^2(\xi_i))c_2^2}{c_1 n}\right)\|\mathbf{e}^x_{k}\|^2\\
			&\quad+\alpha_{k}\frac{10(L^2+\max_{i\in\mathcal{V}}L_i^2(\xi_i))}{c_1 n}\left\|\mathbf{y}_k-\bar{\mathbf{y}}_k\right\|^2\\
			&\le \left(1-\check{\delta}\alpha_k\right)\left(\left(f(\bar{x}_k)-f^*\right)+\left\|\mathbf{e}_k^q+m_k\right\|_2^2\right)+16\alpha_k\left\|m_k\right\|^2+\tilde{P}_k,
		\end{align*}
		}where $\check{\delta}:=\min\left\{\mu,c_2-16\right\}$,
		{\small\begin{align*}
			\tilde{P}_k&:=\alpha_k\left(\frac{8L^2}{n}+\frac{10(L^2+\max_{i\in\mathcal{V}}L_i^2(\xi_i))c_2^2}{c_1 n}\left\|\mathbf{W}-\mathbf{I}\right\|^2\right)\left\|\mathbf{x}_k-\bar{\mathbf{x}}_k\right\|^2+\alpha_{k}\frac{10(L^2+\max_{i\in\mathcal{V}}L_i^2(\xi_i))}{c_1 n}\left\|\mathbf{y}_k-\bar{\mathbf{y}}_k\right\|^2\\
			&\quad+\alpha_k\left(\frac{10\tilde{\gamma}^2}{n}+\frac{10(L^2+\max_{i\in\mathcal{V}}L_i^2(\xi_i))}{c_1 n}\right)\|\mathbf{e}^s_{k}\|^2+\alpha_k\left(\frac{10c_1^2}{n}+\frac{10(L^2+\max_{i\in\mathcal{V}}L_i^2(\xi_i))c_2^2}{c_1 n}\right)\|\mathbf{e}^x_{k}\|^2,
		\end{align*}}the second inequality follows from the condition $c_2>\max\{10(L^2+\max_{i\in\mathcal{V}}L_i^2(\xi_i)),16\}$.
		Then, recursively,
		\begin{align}\label{bar-sum-bound}
			f(\bar{x}_{k+1})-f^*+\left\|\mathbf{e}_{k+1}^q+m_{k+1}\right\|^2
			&\le  \Pi_{t=1}^{k}(1-\check{\delta}\alpha_t)\left(\left(f(\bar{x}_{1})-f^*\right)+\left\|\mathbf{e}_1^q+m_1\right\|^2\right)\notag\\
			&\quad+\sum_{t= 1}^{k}\left(\Pi_{l=t+1}^{k}(1-\check{\delta}\alpha_t)\right)\left(16\alpha_t\|m_t\|^2+\tilde{P}_t\right).
		\end{align}
		The first term on the right hand side of above inequality converges exponentially fast to 0.  Noting that $\epsilon_t$ is a martingale difference sequence and $\lambda_k=c_2\alpha_k$, then by a similar analysis of Proposition \ref{prop:vra-1},
		$$\|m_t\|^2=\mathcal{O}\left(\alpha_t\ln(t+1)\right)~~~\text{a.s.}$$
		In addition, in light of that
		$\left\|\mathbf{e}^s_{t}\right\|^2\le \tilde{c}_s\beta_t\ln (t+1)$, $\left\|\mathbf{e}^x_{t}\right\|^2\le \tilde{c}_x\beta_t\ln (t+1)$,  $\left\|\mathbf{x}_{t+1}-\bar{\mathbf{x}}_{t+1}\right\|^2=\mathcal{O}\left(\alpha_t\ln(t+1)\right)$ and  $\left\|\mathbf{y}_{t+1}-\bar{\mathbf{y}}_{t+1}\right\|^2=\mathcal{O}\left(\alpha_t\ln(t+1)\right)$, we have $$\tilde{P}_t=\mathcal{O}\left(\alpha_t\ln(t+1)\right)~~~\text{a.s.}$$
		Hence, by  \cite[Lemma I2]{Cardot2017fast},
		the last term on the right hand side of (\ref{bar-sum-bound}) can be bounded as
		$$ \sum_{t= 1}^{k}\left(\Pi_{l=t+1}^{k}(1-\check{\delta}\alpha_l)\right)\left(16\alpha_t\|\tilde{m}_t\|^2+\tilde{P}_t\right)=\mathcal{O}\left(\alpha_k\ln(k+1)\right)~~~\text{a.s.}$$
		In summary, 
		$$	f(\bar{x}_{k+1})-f^*+\left\|\mathbf{e}_{k+1}^q+m_{k+1}\right\|^2=\mathcal{O}\left(\alpha_k\ln(k+1)\right)~~~\text{a.s.}$$ 
		Noting that $\|m_{k+1}\|^2=\mathcal{O}\left(\alpha_{k+1}\ln(k+2)\right)$ and $\lim_{k\rightarrow\infty}\frac{\alpha_k\ln(k+1)}{\alpha_{k+1}\ln(k+2)}<\infty$, we have
		\begin{equation*}
			f(\bar{x}_{k+1})-f^*+\left\|\mathbf{e}_{k+1}^q\right\|^2\le 2\left(f(\bar{x}_{k+1})-f^*+\left\|\mathbf{e}_{k+1}^q+m_{k+1}\right\|^2\right)+2\left\|m_{k+1}\right\|^2=\mathcal{O}\left(\alpha_k\ln(k+1)\right).
		\end{equation*}
		Thus, $\mathcal{L}^S_{k+1}=\mathcal{O}\left(\alpha_k\ln(k+1)\right)$ almost surely. The proof is complete.\end{proof}
Theorem \ref{thm:as-DSGT} establishes an almost sure convergence rate of $\mathcal{O}\left(\alpha_k\ln(k+1)\right)$ for VRA-DSGT under diminishing stepsizes, which  may complement the almost sure convergence result on distributed optimization algorithms over the imperfect communication networks \cite{wangGT2022,zhao2023VRA,Zhang2019Sign,Doan2021Quanti}.

The following theorem establishes the convergence rate of VRA-DGT in the almost sure sense under both constant and diminishing stepsizes.

\begin{thm}\label{thm:constant-as}
			Suppose Assumptions \ref{ass:function}, \ref{ass:matrix} and \ref{ass:noise} (with $p_1>2$) hold.
			Then for VRA-DGT, 
			\begin{itemize}
				\item[(i)]  
				$$\mathcal{L}^D_{k+1}=\mathcal{O}\left(\alpha_k\ln(k+1)\right)~~~\text{a.s.}$$
				if $\alpha_k\asymp\frac{1}{k^{b}}$ is nonincreasing, $\gamma_k = c_1 \alpha_k$, $\beta_k \le c_2 \alpha_k$ and that
				 the parameters $b$, $\tilde{\gamma}$, $c_1$ and $c_2$ satisfy $b\in(1/2,1)$
				and the setting (\ref{cond-0}); 
				\item[(ii)]	\begin{align}\label{L-sum-2}
					\mathcal{L}^D_{k+1}
					\le&\left(1-\frac{\delta^D\tilde{\alpha}}{2}\right)^k\mathcal{L}^D_1+\frac{}{}\frac{4\left(\frac{12\tilde{\gamma}}{\eta_w}+\tau^D\right)}{\delta^D}\beta_{k}\ln (k+1)~~~\text{a.s.}
				\end{align}
				if  $\alpha_k\equiv \tilde{\alpha}$, $\gamma_k\equiv c_1\tilde{\alpha}$, $\beta_k\asymp \frac{1}{k^b}$ or $\beta_k=\frac{a_2}{k+a_1}$,  and that the  parameters $b$, $a_1$, $a_2$, $c_1$, $\tilde{\gamma}$ and  $ \tilde{\alpha}$ satisfy the settings (\ref{cond-0}) and (\ref{para-2}), where $\delta^D$ and $\tau^D$ are defined in (\ref{delta-1}) and (\ref{tau-1}). 
\end{itemize}
\end{thm}
\begin{proof}
\textbf{Part (i)}.	By (\ref{y-bound-1}) and (\ref{x-bound}) in Appendix A,
	\begin{align}
		&\left\|\mathbf{y}_{k+1}-\bar{\mathbf{y}}_{k+1}\right\|^2\notag\\
		&\le (1-\tilde{\gamma} \eta_w)\left\|\mathbf{y}_{k}-\bar{\mathbf{y}}_{k}\right\|^2+\frac{2L^2}{\tilde{\gamma} \eta_w}\left\|\mathbf{x}_{k+1}-\mathbf{x}_{k}\right\|^2+\frac{4\tilde{\gamma}}{ \eta_w}\left\|\mathbf{e}_{k}^s\right\|^2+\frac{4\tilde{\gamma}}{ \eta_w}\left\|\mathbf{e}_{k+1}^s\right\|^2\notag\\
		&\le \left(1-\tilde{\gamma} \eta_w+\frac{10L^2}{\tilde{\gamma} \eta_w}\alpha_k^2\right)\left\|\mathbf{y}_{k}-\bar{\mathbf{y}}_{k}\right\|^2+\frac{10L^2}{\tilde{\gamma} \eta_w}\left(\gamma_k^2\left\|\mathbf{e}_k^x\right\|^2+\gamma_k^2\left\|\mathbf{W}-\mathbf{I}\right\|^2\left\|\mathbf{x}_k-\bar{\mathbf{x}}_k\right\|^2\right.\notag\\
		&\quad\left.+n\alpha_k^2\left\|\bar{q}_k\right\|_2^2\right)+\frac{2\tilde{\gamma}\left(2+5L^2\alpha_k^2\right)}{\eta_w}\left\|\mathbf{e}_{k}^s\right\|^2+\frac{4\tilde{\gamma}}{ \eta_w}\left\|\mathbf{e}_{k+1}^s\right\|^2.\label{y-bound-2}
	\end{align}
	The above inequality combining with the upper bounds of  $\left\|\mathbf{x}_{k+1}-\bar{\mathbf{x}}_{k+1}\right\|^2$ and $f(\bar{x}_{k+1})-f^*$ in Lemmas \ref{lem:xstar-1} and \ref{lem:x-1} (regardless of the conditional expectation) implies 
	\begin{align}\label{L-bound-GT-1}
		\mathcal{L}^D_{k+1}
		\le&(1-\mu\alpha_k)(f(\bar{x}_k)-f^*)+\left[1-\left(\frac{\tilde{\gamma} \eta_w}{\alpha_k}-\frac{2}{c_1 \eta_w}\right)\alpha_k+\left(\frac{10L^2}{\tilde{\gamma} \eta_w}+30L^2\right)\alpha_k^2\right]\left\|\mathbf{y}_{k}-\bar{\mathbf{y}}_{k}\right\|^2\notag\\
		&+\left[1-\left(c_1 \eta_w-\frac{8L^2}{n}\right)\alpha_k+\left(\frac{10L^2}{ \tilde{\gamma}\eta_w}+30L^2\right)\left\|\mathbf{W}-\mathbf{I}\right\|^2c_1^2\alpha_k^2\right]\left\|\mathbf{x}_k-\bar{\mathbf{x}}_k\right\|^2\notag\\
		&+\left[\left(\frac{10L^2}{ \tilde{\gamma}\eta_w}+30L^2\right)c_1^2+\frac{2c_1}{ \eta_w}+\frac{10c_1^2}{n}\right]\alpha_k\left\|\mathbf{e}^x_{k}\right\|^2+\left(\frac{4\tilde{\gamma}+10\tilde{\gamma}\alpha_kL^2}{\eta_w}+\frac{10\tilde{\gamma}^2}{n}\right)\alpha_k\left\|\mathbf{e}^s_{k}\right\|^2\notag\\
		&+\frac{4\tilde{\gamma}}{ \eta_w}\left\|\mathbf{e}_{k+1}^s\right\|^2\notag\\
		\le&(1-\mu\alpha_k)(f(\bar{x}_k)-f^*)+\left[1-\left(\frac{\tilde{\gamma} \eta_w}{\alpha_k}-\frac{2}{c_1 \eta_w}\right)\alpha_k+\left(\frac{10L^2}{\tilde{\gamma} \eta_w}+30L^2\right)\alpha_k^2\right]\left\|\mathbf{y}_{k}-\bar{\mathbf{y}}_{k}\right\|^2\notag\\
		&+\left[1-\left(c_1 \eta_w-\frac{8L^2}{n}\right)\alpha_k+\left(\frac{10L^2}{ \tilde{\gamma}\eta_w}+30L^2\right)\left\|\mathbf{W}-\mathbf{I}\right\|^2c_1^2\alpha_k^2\right]\left\|\mathbf{x}_k-\bar{\mathbf{x}}_k\right\|^2\notag\\
		&+\left[\left(\frac{10L^2}{ \tilde{\gamma}\eta_w}+30L^2\right)c_1^2+\frac{2c_1}{ \eta_w}+\frac{10c_1^2}{n}\right]\alpha_k\tilde{c}_x\beta_k\ln (k+1)\notag\\
		&+\left(\frac{4\tilde{\gamma}+10\tilde{\gamma}\alpha_kL^2}{\eta_w}+\frac{10\tilde{\gamma}^2}{n}\right)\alpha_k\tilde{c}_s\beta_k\ln (k+1)+\frac{4\tilde{\gamma}}{ \eta_w}\tilde{c}_s\beta_{k+1}\ln (k+2)~~~\text{a.s.}
		%		\notag\\
		%		\le&\left(1-\delta^D\alpha_k+\tilde{\delta}^D\alpha_k^2\right)\mathbb{E}\left[\mathcal{L}^D_k\right]+\left(\frac{12\tilde{\gamma}}{\eta_w}+\tau^D\right)\alpha_k\beta_{k}\notag\\
		%		\le&\left(1-\frac{\delta^D\alpha_k}{2}\right)\mathbb{E}\left[\mathcal{L}^D_k\right]+\left(\frac{12\tilde{\gamma}}{\eta_w}+\tau^D\right)\alpha_k\beta_{k},		
	\end{align}
	where the second inequality follows from Proposition \ref{prop:vra-1}.
	In light of $\beta_k$  being nonincreasing and $\ln (k+2)\le 2\ln (k+1)$, the last term on the right hand side of (\ref{L-bound-GT-1}) can be bounded as
	\begin{align*}
		\frac{4\tilde{\gamma}}{ \eta_w}\tilde{c}_s\beta_{k+1}\ln (k+2)=\frac{4\tilde{\gamma}}{\eta_w}\tilde{c}_s\beta_k\ln (k+1)\frac{\beta_{k+1}\ln (k+2)}{\beta_k\ln (k+1)}\le \frac{8\tilde{\gamma}}{\eta_w}\tilde{c}_s\beta_k\ln (k+1).
	\end{align*}
	Subsequently, (\ref{L-bound-GT-1}) reduces to
	\begin{align}\label{L-bound-1}
		\mathcal{L}^D_{k+1}
		\le&\left(1-\delta^D\alpha_k+\tilde{\delta}^D\alpha_k^2\right)\mathcal{L}^D_k+\left(\frac{12\tilde{\gamma}}{\eta_w}+\tau^D\right)\alpha_k\beta_{k}\ln (k+1)	~~~\text{a.s.}	
	\end{align}
	where $\delta^D$, $\tilde{\delta}^D$ and $\tau^D$ are defined in  (\ref{delta-1}) and (\ref{tau-1}). Then there exists a large enough integer $k_0$ such that
	\begin{align*}
		\mathcal{L}^D_{k+1}
		\le&\left(1-\frac{\delta^D}{2}\alpha_k\right)\mathcal{L}^D_{k}+\left(\frac{12\tilde{\gamma}}{\eta_w}+\tau^D\right)c_3\alpha_k^2\ln (k+1)	~~~\text{a.s.}	,~\forall k\ge k_0.
	\end{align*}
	Recursively, we have
	{\small\begin{align*}
			\mathcal{L}^D_{k+1}
			\le&\Pi_{t=k_0}^k\left(1-\frac{\delta^D}{2}\alpha_t\right)\mathcal{L}^D_1+\left(\frac{12\tilde{\gamma}}{\eta_w}+\tau^D\right)c_3\ln (k+1)\sum_{t= k_0}^k\Pi_{l=t+1}^k\left(1-\frac{\delta^D}{2}\alpha_l\right)\alpha_t^2	~~~\text{a.s.}, ~\forall k\ge k_0.
	\end{align*}}
The above inequality, 	together with (\ref{L-sum-0}) and \cite[Lemma I2]{Cardot2017fast}, implies  $$\mathcal{L}^D_{k+1}=\mathcal{O}\left(\alpha_k\ln(k+1)\right)~~~\text{a.s.}$$
\textbf{Part (ii)}.	Given  fact $\alpha_k \equiv \tilde{\alpha}$, $\tilde{\alpha} \le \frac{2\tilde{\delta}^D}{\delta^D}$,  (\ref{L-bound-1}) reduces to
	\begin{align*}
		\mathcal{L}^D_{k+1}
		\le&\left(1-\frac{\delta^D\tilde{\alpha}}{2}\right)\mathcal{L}^D_k+\left(\frac{12\tilde{\gamma}}{\eta_w}+\tau^D\right)\tilde{\alpha}\beta_{k}\ln (k+1)~~~\text{a.s.}		
	\end{align*}
	Recursively, 
	\begin{align*}
		\mathcal{L}^D_{k+1}
		\le&\left(1-\frac{\delta^D\tilde{\alpha}}{2}\right)^k\mathcal{L}^D_1+\left(\frac{12\tilde{\gamma}}{\eta_w}+\tau^D\right)\tilde{\alpha}\sum_{t=1}^k\left(1-\frac{\delta^D\tilde{\alpha}}{2}\right)^{k-t}\beta_{t}\ln(t+1)~~~\text{a.s.}
	\end{align*}
	 By condition $\sup_{k\ge 1}\frac{\beta_{k}^2}{\beta_{k+1}^2}\le 1+\frac{\delta^D\tilde{\alpha}}{4-2\delta^D\tilde{\alpha}}$ and \cite[Lemma 7]{Li2022Perform}, the last term on the right hand side of above inequality has $$\sum_{t=1}^k\left(1-\frac{\delta^D\tilde{\alpha}}{2}\right)^{k-t}\beta_{t}\ln(t+1)\le\ln(k+1)\sum_{t=1}^k\left(1-\frac{\delta^D\tilde{\alpha}}{2}\right)^{k-t}\beta_{t}\le \frac{4}{\delta^D\tilde{\alpha}}\beta_{k}\ln(k+1).$$
	  Then, we arrive at (\ref{L-sum-2}). The proof is complete.
\end{proof}

	 Theorem \ref{thm:constant-as} establishes both asymptotic and non-asymptotic convergence rates for VRA-DGT in the almost sure sense under diminishing and constant stepsizes. By choosing $\beta_k\asymp\frac{1}{k^{b}}$, $\alpha_k\asymp\frac{1}{k^{b}}$  or $\alpha_k\equiv\tilde{\alpha}$,   VRA-DGT achieves the convergence rate of $\mathcal{O}\left(\frac{\ln (k+1)}{k^{b}}\right)$ in the almost sure sense. To our best
		knowledge, this is the first result that establishes an almost sure convergence rate for distributed
		optimization algorithms under both imperfect communication and constant stepsizes.

\section{Experimental Results}\label{sec:num-experi}

We evaluate the empirical performance of the VRA-DGT and VRA-DSGT algorithms on a logistic regression problem \cite{Iakovi2023S-NEAR} defined as follows:
\begin{equation*}
	\min_{x\in\mathbb{R}^d} f(x)\define\frac{1}{N}\sum_{s=1}^N \log\left(1+\e^{-b_s\langle a_s,x\rangle}\right)+\psi\|x\|^2,
\end{equation*}
where $\psi=\frac{1}{N}$ is the regularization parameter and $\{(a_1,b_1)$, $\cdots$, $(a_N,b_N)\}$ is the dataset. We employ a group of $n$ agents to collaboratively address the logistic regression problem. Using the Dirichlet distribution, each agent $j \in \mathcal{V}$ is assigned a distinct non-iid subset $\mathcal{S}_j$ from the dataset. Consequently, the private function for each agent $j\in\mathcal{V}$ is given by:
\begin{equation*}
f_j(x)= \frac{1}{|\mathcal{S}_j|}\sum_{s\in \mathcal{S}_j} \log\left(1+\e^{-b_s\langle a_s,x\rangle}\right)+\psi\|x\|^2.
\end{equation*}
For $\forall i\in \mathcal{V}$, the $(i,j)$-th entry of the weight matrix 
 $\mathbf{W}$ is set to 
\begin{equation*}
	w_{ij}=\left\{
	\begin{aligned}
		&\frac{1}{\max\left\{|\mathcal{N}_{i}|,|\mathcal{N}_{j}|\right\}+1},\quad j\in \mathcal{N}_{i},\\
		&1-\sum_{j\in\mathcal{N}_{i}}w_{ij},\quad j=i,
	\end{aligned}\right.
\end{equation*}
where $|\mathcal{N}_{i}|$ is the  cardinality of $\mathcal{N}_{i}$. There are two inexact communication scenarios considered, including
\begin{itemize}
 	\item [1.]\textbf{Noisy channel \cite{Reisizadeh2023Vary}}. Assume communication between agents occurs over a Gaussian channel. Specifically, when
	agent $j$ sends its state $x_{j,k}$ to its neighbors, the received signal at its neighbors is $x_{j,k}+\zeta_{j,k}$, where $\zeta_{j,k}$ is a
	zero-mean Gaussian noise with variance $\sigma_\zeta^2$, which is independent across $(i,j)$ and $k$.
	In our experiments, the variance $\sigma^2_\zeta$ is set to $\sigma^2_\zeta=1,~50~\text{and}~100$.
	\item[2.] \textbf{Probabilistic quantizer \cite{YUAN2012DDA}}.  Assume the state of an agent is quantized to a specific number of bits using a probabilistic quantizer before transmission to its neighbors. Specifically,  for $\theta\in\mathbb{R}$, the quantized value $Q_{\Delta}(\theta)$ is given by
	\begin{equation*}
		Q_{\Delta}(\theta)=\left\{
		\begin{aligned}
			& \lfloor \theta \rfloor \quad \text{with probability}(\lceil \theta \rceil-\theta)\Delta,\\
			&\lceil \theta \rceil\quad \text{with probability}(\theta-\lfloor \theta \rfloor)\Delta,
		\end{aligned}\right.
	\end{equation*}
	where parameter $\Delta$ is a positive real value, $\lfloor \theta \rfloor$  and $\lceil \theta \rceil$ denote the operations of rounding down and up to the nearest integer multiple of $\frac{1}{\Delta}$, respectively. For any $x\in\mathbb{R}^d$, the resulting error term, denoted as $\xi\define x-Q_{\Delta}(x)$, possesses a zero mean and bounded variance, that is, $\mathbb{E}\left[\xi\right]=0$ and $\mathbb{E}\left[\left\|\xi\right\|^2\right]\le \frac{d}{4\Delta^2}$ \cite{Iakovi2023S-NEAR}. In our experiments, the parameter $\Delta$ is set to $\Delta=1,~10~\text{and}~100$.
\end{itemize}
{\small\begin{table}[http]\label{tab:parameters-1}
		\renewcommand{\tablename}
		\caption{\centering{ Table 2. Stepsize and parameter settings for Mushrooms dataset.} \protect \\  }
		\centering
		\begin{tabular}{c|c|c|c|c}
			\hline
			Methods& Stepsizes&\multicolumn{3}{c}{Parameters}\\
			\hline
			Algorithm in \cite{wangGT2022}&$\frac{0.05}{1+0.001\times k^{0.8}}$&$\gamma_k=\frac{10^{-4}}{1+0.001\times k^{0.6}}$&-&-\cr
			\hline
			VRA-GT&$\frac{0.05}{1+0.001\times k^{0.8}}$&$\gamma_k=\frac{10^{-4}}{1+0.001\times k^{0.6}}$&$\beta_{k}=\frac{1\times 10^{-7}}{1+k}$&$\tilde{\gamma}=0.99$\cr
			\hline
			VRA-DGT&$\frac{0.05}{1+0.001\times k^{0.8}}$&$\gamma_k=\frac{1}{1+0.001\times k^{0.8}}$&$\beta_k=\frac{1\times 10^{-7}}{1+k}$&$\tilde{\gamma}=0.99$\cr
			\hline
			IC-GT&0.05&$\tilde{\gamma}=10^{-7}$&-&-\cr
			\hline
			VRA-DSGT&$\frac{0.05}{1+0.001\times k^{0.8}}$&$\gamma_k=\frac{1}{1+0.001\times k^{0.8}}$& {\makecell[c]{$\beta_k=\frac{1\times 10^{-7}}{1+k}$, \\ $\lambda_k=\frac{1}{1+0.001\times k^{0.8}}$}}  &$\tilde{\gamma}=0.99$\cr
			\hline
		\end{tabular}
\end{table}}
\begin{figure*}[htb]
	\centering
	\subfigure[Noisy channel.]{
		\includegraphics[width=2.8in]{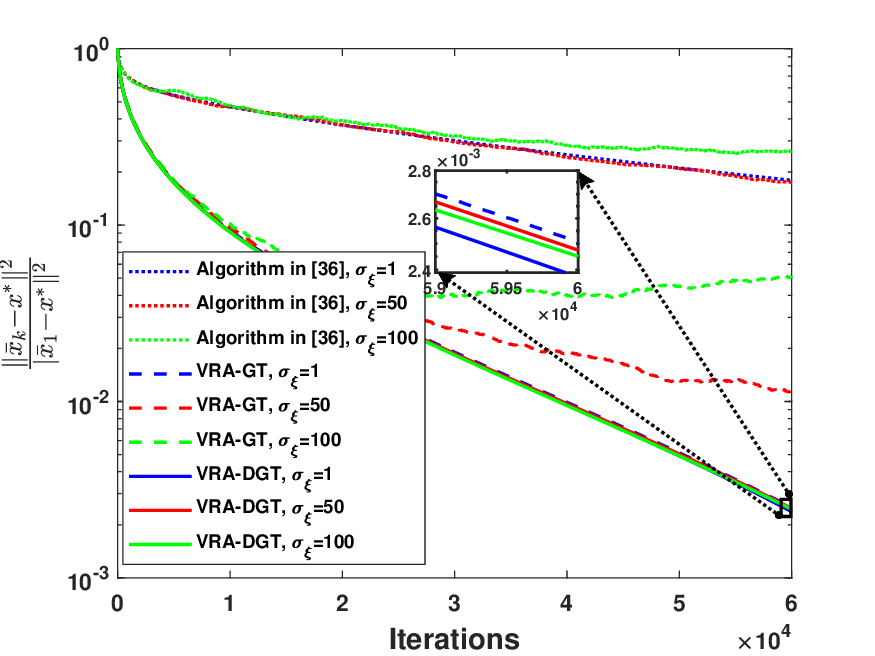}
	}\hspace{-1mm}
	\subfigure[Probabilistic quantizer.]{
		\includegraphics[width=2.8in]{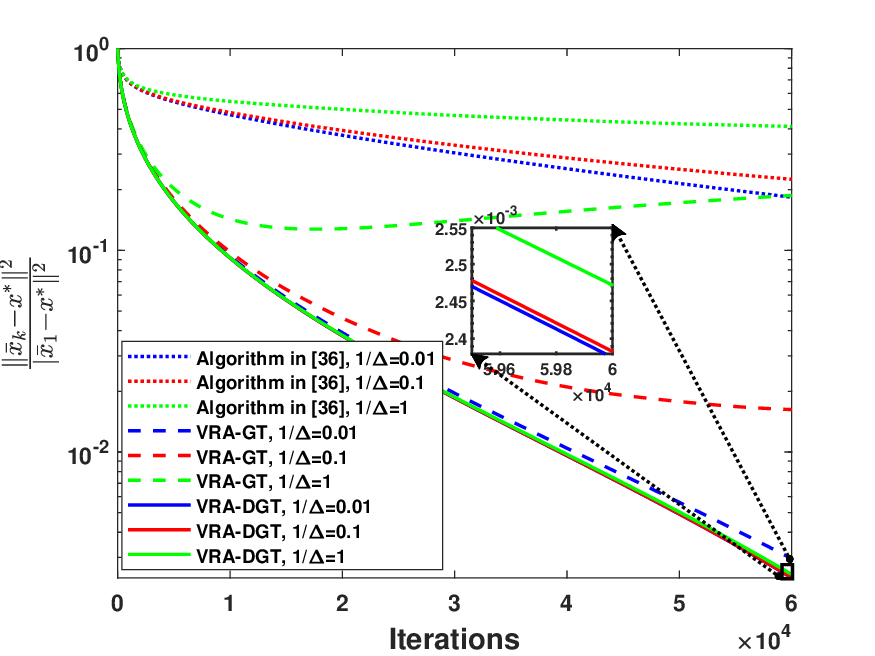}
	}
	\caption{{\small  The value of  $\frac{\|\bar{x}_k-x^*\|^2}{\|\bar{x}_1-x^*\|^2}$ w.r.t  the number of iterations (VRA-DGT).}}
	\label{fig-1}
\end{figure*}

\begin{figure*}[htb]
	\centering
	\subfigure[Noisy channel.]{
		\includegraphics[width=2.8in]{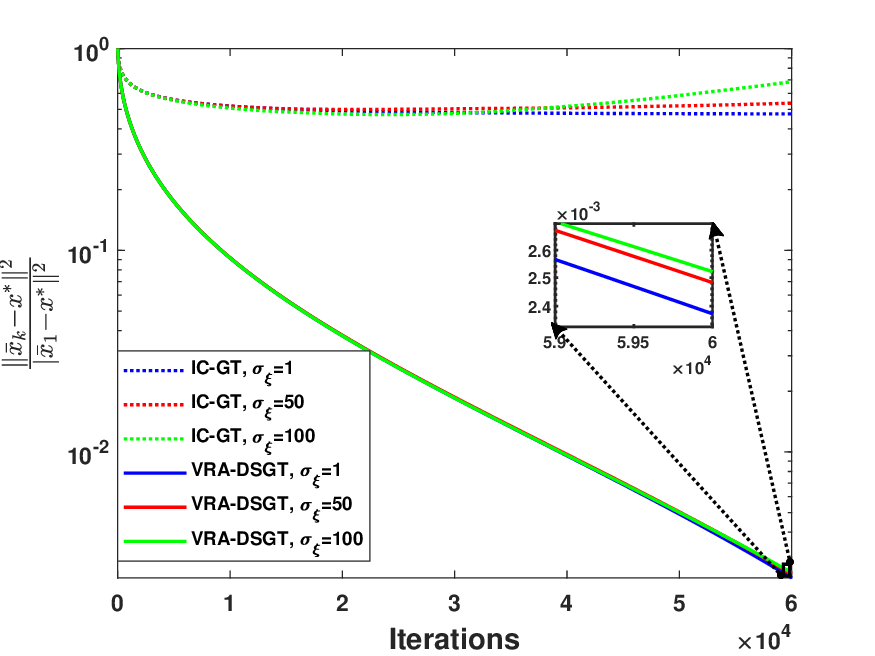}
	}\hspace{-1mm}
	\subfigure[Probabilistic quantizer.]{
		\includegraphics[width=2.8in]{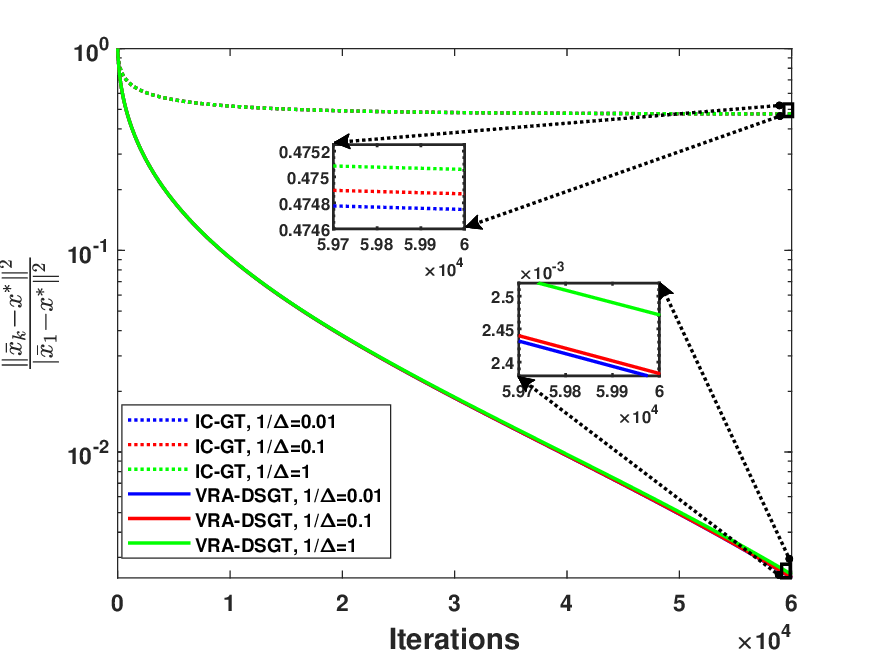}
	}
	\caption{{\small  The value of  $\frac{\|\bar{x}_k-x^*\|^2}{\|\bar{x}_1-x^*\|^2}$ w.r.t  the number of iterations (VRA-DSGT).}}
	\label{fig-2}
\end{figure*}

We compare the performance of   VRA-DGT and VRA-DSGT with the algorithms including VRA-GT \cite{zhao2023VRA}, the algorithm proposed in \cite{wangGT2022} and IC-GT \cite{shah2023SGT} under  deterministic and stochastic gradient conditions.  Their convergence performance is evaluated on the Mushrooms \cite{dua2017uci} and Adult \cite{adult_2} datasets respectively using $\frac{\|\bar{x}_k-x^*\|^2}{\|\bar{x}_1-x^*\|^2}$.

\textbf{Mushrooms dataset \cite{dua2017uci}, random graph.} The graph $\mathcal{G}$ is generated by adding random links to a 50-agent ring network with probability $p=0.3$. The algorithms' stepsizes and parameters are detailed in Table 2. The convergence performance of   VRA-DGT, VRA-GT, and the algorithm proposed in \cite{wangGT2022} for distributed deterministic optimization problems is presented in Figure \ref{fig-1}.
Figure \ref{fig-1} illustrates that VRA-DGT, VRA-GT, and the algorithm proposed in \cite{wangGT2022} achieve convergence under varying levels of channel noise and communication quantization. The algorithms that incorporate the VRA mechanism—VRA-GT and VRA-DGT—exhibit more effective convergence than the one presented in \cite{wangGT2022}. Moreover, VRA-DGT exhibits more stable convergence as the noise level increases compared to VRA-GT. This enhanced stability is likely due to VRA-DGT applying the VRA mechanism to both decision variable updating and cumulative global gradient tracking, while VRA-GT incorporates the VRA mechanism only in cumulative global gradient tracking.
The convergence performance of VR-DSGT and  IC-GT for distributed stochastic optimization problems is presented in Figure \ref{fig-2}. Figure \ref{fig-2} shows that both IC-GT and VRA-DSGT achieve  convergence across different levels of channel noise and communication quantization. Moreover, VRA-DSGT demonstrates better convergence performance compared to IC-GT.

\begin{figure*}[htb]
	\centering
	\subfigure[Noisy channel.]{
		\includegraphics[width=2.8in]{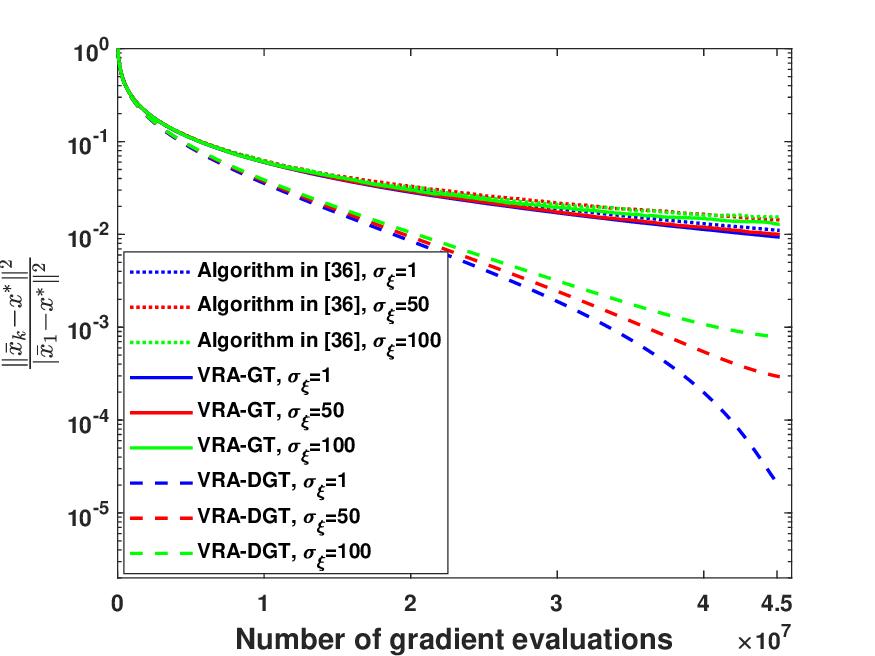}
	}\hspace{-1mm}
	\subfigure[Probabilistic quantizer.]{
		\includegraphics[width=2.8in]{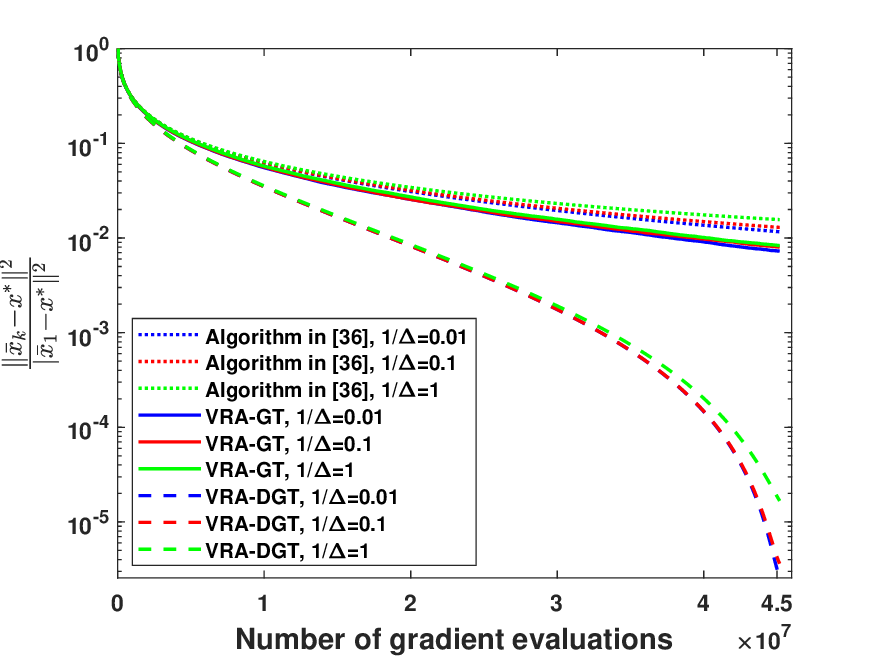}
	}
	\caption{{\small  The value of  $\frac{\|\bar{x}_k-x^*\|^2}{\|\bar{x}_1-x^*\|^2}$ w.r.t  the number of iterations (VRA-DGT).}}
	\label{fig-3}
\end{figure*}

\begin{figure*}[htb]
	\centering
	\subfigure[Noisy channel.]{
		\includegraphics[width=2.8in]{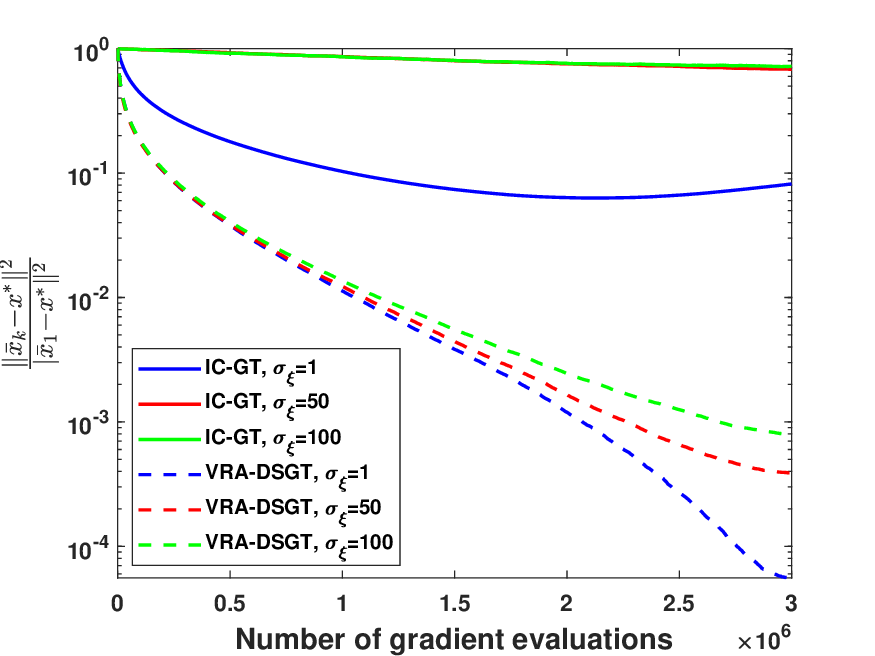}
	}\hspace{-1mm}
	\subfigure[Probabilistic quantizer.]{
		\includegraphics[width=2.8in]{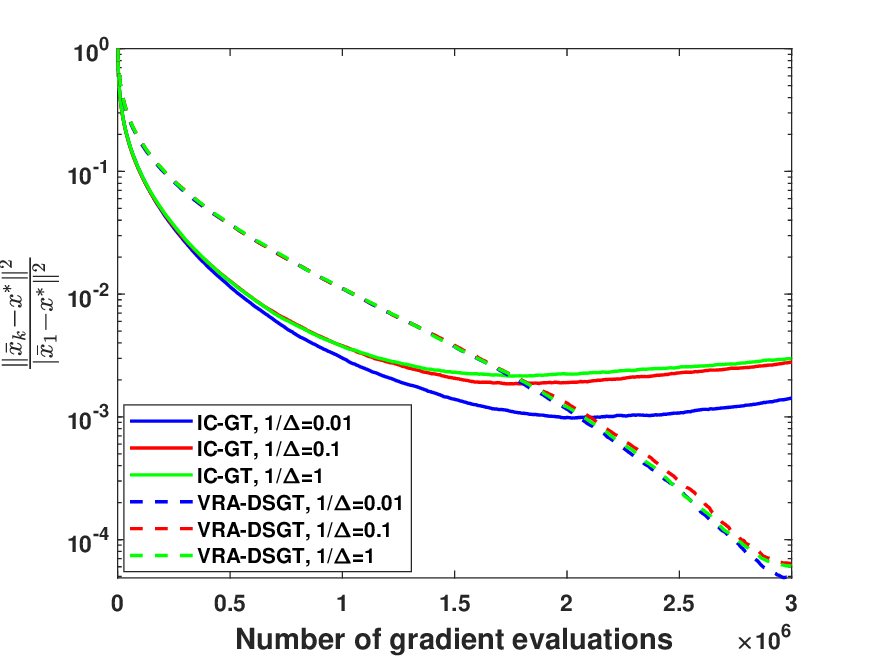}
	}
	\caption{{\small  The value of  $\frac{\|\bar{x}_k-x^*\|^2}{\|\bar{x}_1-x^*\|^2}$ w.r.t  the number of gradient evaluations (VRA-DSGT).}}
	\label{fig-4}
\end{figure*}

\textbf{Adult dataset  \cite{adult_2}, grid graph.}
{\small\begin{table}[http]\label{tab:parameters-2}
		\renewcommand{\tablename}
		\caption{\centering{ Table 3. Stepsize and parameter settings for Adult dataset.} \protect \\  }
		\centering
		\begin{tabular}{c|c|c|c|c}
			\hline
			Methods& Stepsizes&\multicolumn{3}{c}{Parameters}\\
			\hline
			Algorithm in \cite{wangGT2022}&$\frac{\tilde{\alpha}}{1+0.01\times k^{0.8}}$&$\gamma_k=\frac{\hat{\gamma}}{1+0.01\times k^{0.6}}$&-&-\cr
			\hline
			VRA-GT&$\frac{\tilde{\alpha}}{1+0.01\times k^{0.8}}$&$\gamma_k=\frac{\hat{\gamma}}{1+0.01\times k^{0.6}}$&$\beta_{k}=\frac{10^{-4}}{k^{0.9}}$&$\tilde{\gamma}=0.01$\cr
			\hline
			VRA-DGT&$\tilde{\alpha}$&$\hat{\gamma}$&$\beta_k=\frac{10^{-4}}{k^{0.9}}$&$\tilde{\gamma}=0.01$\cr
			\hline
			IC-GT&$\tilde{\alpha}$&$\hat{\gamma}$&-&-\cr
			\hline
			VRA-DSGT&$\frac{\tilde{\alpha}}{1+0.01\times k^{0.6}}$&$\gamma_k=\frac{\hat{\gamma}}{1+0.01\times k^{0.6}}$& {\makecell[c]{$\beta_k=\frac{10^{-4}}{k^{0.9}}$, \\ $\lambda_k=\frac{1}{1+0.01\times k^{0.6}}$}}  &$\tilde{\gamma}=0.01$\cr
			\hline
		\end{tabular}
\end{table}}
 The graph $\mathcal{G}$ is constructed as a $10 \times 10$ grid with 100 agents. The algorithms' stepsizes and parameters are detailed in Table 3. Here, $\beta_k$, $\tilde{\gamma}$, and $\lambda_k$ are fixed to the same values for VRA-D(S)GT and VRA-GT. The  step size $\alpha_k$ and parameter $\gamma_k$ for all algorithms are tuned by first selecting their magnitudes according to the original papers, and then optimizing the constants $\tilde{\alpha}$ and $\tilde{\gamma}$ from $\{10^{-4},10^{-3},10^{-2},10^{-1},1\}$.
The convergence performance of VRA-DGT, VRA-GT, and the algorithm in \cite{wangGT2022} is presented in Figure \ref{fig-3}, while that of VRA-DSGT and IC-GT is presented in Figure \ref{fig-4}, where, in both figures, performance is plotted as loss versus gradient evaluations. Figure \ref{fig-3} illustrates that VRA-DGT maintains the better convergence performance across information-sharing noise types and levels, which also confirms Theorem 2 that VRA-DGT converges under constant stepsize. 
For noisy channels, Figure \ref{fig-4}(a) shows that VRA-DGT performs better than IC-GT under different levels of channel noise. For communication quantization, Figure \ref{fig-4}(b) shows that IC-GT converges faster in the initial stage, while VRA-DSGT achieves higher final accuracy and more stable convergence. The underlying reason may be that IC-GT uses a constant stepsize and a constant noise attenuation parameter.

%\textbf{Conclusion and discussion.}

\textbf{Conclusion.} This paper proposes the  VRA  based distributed (stochastic) gradient tracking algorithm for distributed (stochastic) optimization problems over the imperfect-communication undirected networks. For  strongly convex and smooth  objective functions, VRA-D(S)GT releases the two-time-scale condition on stepsizes and noise attenuation parameters required in existing works. VRA-D(S)GT achieves a convergence rate of $\mathcal{O}(k^{-1})$ in the mean square sense and $\mathcal{O}\left(\frac{\ln(k+1)}{k^b}\right)$, $\forall b \in (0.5, 1)$ in the almost sure sense.   Future work will extend VRA-D(S)GT to more general distributed optimization problems over directed networks. Another important area pertains to uncertainty quantification for VRA-D(S)GT.

\textbf{Acknowledgment.}  The authors are  grateful to the referees of \cite{zhao2023VRA} for  suggesting us to work on this topic. The authors would like to thank two anonymous reviewers for their valuable comments.
The research is supported by National Key R$\&$D Program of China No. 2023YFA1009200, the NSFC \#12401418, the NSFC \#12471283 and  Fundamental Research Funds for the Central Universities DUT24LK001.

\textbf{Data availability}

The authors confirm that all data generated or analyzed during this study are included in this
published article.

\textbf{Declarations}

The authors have no competing interests to declare that are relevant to the content of this article.

\bibliographystyle{siam}
\bibliography{mybib}

\section*{Appendix A: Proofs of Lemmas \ref{lem:vr}-\ref{lem:xstar-1}}

\noindent\textbf{A.1. Proof of Lemma \ref{lem:vr}}
\begin{proof}
	By the definition of $\mathbf{e}^q_k$ and the recursion (\ref{alg:q}),
	\begin{equation}\label{eq-recur}
		\mathbf{e}_{k+1}^q=(1-\lambda_k)\mathbf{e}_{k}^q+(1-\lambda_k)\left(\mathbf{g}_k-\tilde{\mathbf{g}}_k\right)+\tilde{\mathbf{g}}_{k+1}^{'}-\mathbf{g}_{k+1}.
	\end{equation}
	Taking square
	norms and conditional expectation on both sides of the above equality, 
	\begin{align}
		&\mathbb{E}\left[\left\|\mathbf{e}_{k+1}^q\right\|^2\big|\mathcal{F}_k^s\right]\notag\\
		=&(1-\lambda_k)^2\left\|\mathbf{e}_{k}^q\right\|^2+2\mathbb{E}\left[\langle(1-\lambda_k)\mathbf{e}_{k}^q, (1-\lambda_k)\left(\mathbf{g}_k-\tilde{\mathbf{g}}_k\right)+\tilde{\mathbf{g}}_{k+1}^{'}-\mathbf{g}_{k+1}\rangle\big|\mathcal{F}_k^s\right]\notag\\
		&+\mathbb{E}\left[\left\|(1-\lambda_k)\left(\mathbf{g}_k-\tilde{\mathbf{g}}_k\right)+\tilde{\mathbf{g}}_{k+1}^{'}-\mathbf{g}_{k+1}\right\|^2\big|\mathcal{F}_k^s\right]\notag\\
		=&(1-\lambda_k)^2\left\|\mathbf{e}_{k}^q\right\|^2+\mathbb{E}\left[\left\|(1-\lambda_k)\left(\mathbf{g}_k-\tilde{\mathbf{g}}_k\right)+\tilde{\mathbf{g}}_{k+1}^{'}-\mathbf{g}_{k+1}\right\|^2\big|\mathcal{F}_k^s\right]\notag\\
		\le &(1-\lambda_k)^2\left\|\mathbf{e}_{k}^q\right\|^2+3(1-\lambda_k)^2\mathbb{E}\left[\left\|\tilde{\mathbf{g}}_{k+1}^{'}-\tilde{\mathbf{g}}_k\right\|^2\big|\mathcal{F}_k^s\right]+3(1-\lambda_k)^2\mathbb{E}\left[\left\|\mathbf{g}_{k+1}-\mathbf{g}_k\right\|^2\big|\mathcal{F}_k^s\right]\notag\\
		&+3\lambda_k^2\mathbb{E}\left[\left\|\tilde{\mathbf{g}}_{k+1}^{'}-\mathbf{g}_{k+1}\right\|^2\big|\mathcal{F}_k^s\right]\notag\\
		\le &(1-\lambda_k)^2\left\|\mathbf{e}_{k}^q\right\|^2+6(1-\lambda_k)^2L^2\mathbb{E}\left[\left\|\mathbf{x}_{k+1}-\mathbf{x}_k\right\|^2\big|\mathcal{F}_k^s\right]+3\lambda_k^2n\sigma_\xi^2,\label{eq-bound}
	\end{align}
	where the second equality follows from $\mathbb{E}\left[\tilde{\mathbf{g}}_k\big|\mathcal{F}_k^s\right]=\mathbf{g}_k$ and $\mathbb{E}\left[\tilde{\mathbf{g}}^{'}_{k+1}\big|\mathcal{F}_k^s\right]=\mathbb{E}\left[\mathbb{E}\left[\tilde{\mathbf{g}}^{'}_{k+1}\big|\mathcal{F}_{k+1}^x\right]\big|\mathcal{F}_k^s\right]=\mathbb{E}\left[\mathbf{g}_{k+1}\big|\mathcal{F}_k^s\right]$, the last inequality follows from the Lipschitz smoothness of $ F_i(\cdot;\xi_i)$ and  the condition  $\mathbb{E}\left[\left\|\nabla F_i(x;\xi_i)-\nabla f_i(x)\right\|^2\right]\le  \sigma_\xi^2$. By recursion (\ref{alg:x}) and the fact $\left(\mathbf{W}-\mathbf{I}\right)\bar{\mathbf{x}}_k=\mathbf{0}$,
	\begin{align}\label{x-bound}
		&\mathbb{E}\left[\left\|\mathbf{x}_{k+1}-\mathbf{x}_k\right\|^2\big|\mathcal{F}_k^s\right]\notag\\
		&=\mathbb{E}\left[\left\|\gamma_k\left(\mathbf{W}-\mathbf{I}\right)(\mathbf{x}_k-\bar{\mathbf{x}}_k)+\gamma_k\mathbf{e}_k^x-\alpha_k\mathbf{y}_k\right\|^2\big|\mathcal{F}_k^s\right]\notag\\
		&=\mathbb{E}\left[\left\|\gamma_k\left(\mathbf{W}-\mathbf{I}\right)(\mathbf{x}_k-\bar{\mathbf{x}}_k)+\gamma_k\mathbf{e}_k^x-\alpha_k(\mathbf{y}_k-\bar{\mathbf{y}}_k)-\alpha_k\mathbf{1}(\bar{q}_{k}+\tilde{\gamma}\bar{e}^s_{k})^\intercal\right\|^2\big|\mathcal{F}_k^s\right]\notag\\
		&\le 5\gamma_k^2\left\|\mathbf{W}-\mathbf{I}\right\|^2\left\|\mathbf{x}_k-\bar{\mathbf{x}}_k\right\|^2+5\gamma_k^2\mathbb{E}\left[\left\|\mathbf{e}_k^x\right\|^2\big|\mathcal{F}_k^s\right]+5\alpha_k^2\left\|\mathbf{y}_k-\bar{\mathbf{y}}_k\right\|^2\notag\\
		&\quad+5\alpha_k^2n\left\|\bar{q}_k\right\|_2^2+5\alpha_k^2\tilde{\gamma}^2\left\|\mathbf{e}_k^s\right\|^2,
	\end{align}
	where
	\begin{equation*}
			\bar{q}_{k}:=\frac{1}{n} \sum_{j=1}^nq_{j,k}, ~\bar{e}_{k}^s:=\frac{1}{n} \sum_{j=1}^n\left(z_{j,k}-\sum_{j\in \mathcal{N}_{i}}s_{j,k}\right),
	\end{equation*}
	the second equality follows from the fact $\bar{y}_{k}
	=\bar{q}_{k}+\tilde{\gamma}\bar{e}^s_{k}$. Then, taking above inequality  into (\ref{eq-bound}), we arrive at
	\begin{align*}
		\mathbb{E}\left[\left\|\mathbf{e}_{k+1}^q\right\|^2\big|\mathcal{F}_k^s\right]
		&\le (1-\lambda_k)^2\left\|\mathbf{e}_{k}^q\right\|^2+30L^2(1-\lambda_k)^2\left(\gamma_k^2\left\|\mathbf{W}-\mathbf{I}\right\|^2\left\|\mathbf{x}_k-\bar{\mathbf{x}}_k\right\|^2+\gamma_k^2\mathbb{E}\left[\left\|\mathbf{e}_k^x\right\|^2\big|\mathcal{F}_k^s\right]\right.\notag\\
		&\left.\quad+\alpha_k^2\left\|\mathbf{y}_k-\bar{\mathbf{y}}_k\right\|^2+\alpha_k^2n\left\|\bar{q}_k\right\|_2^2+5\alpha_k^2\tilde{\gamma}^2\left\|\mathbf{e}_k^s\right\|^2\right)+3\lambda_k^2n\sigma_\xi^2.
	\end{align*}
\end{proof}

\noindent\textbf{A.2. Proof of Lemma \ref{lem:y-1}}
\begin{proof}
	Note that $\bar{\mathbf{y}}_{k}=\mathbf{1}\bar{y}_{k}^\intercal=\frac{\mathbf{1}\mathbf{1}^\intercal}{n}\mathbf{y}_{k}$, we have
	\begin{align*}
		\mathbf{y}_{k+1}-\bar{\mathbf{y}}_{k+1}&=\left(\mathbf{W}_{\tilde{\gamma}}-\frac{\mathbf{1}\mathbf{1}^\intercal}{n}\right)\mathbf{y}_{k}+\left(\mathbf{I}-\frac{\mathbf{1}\mathbf{1}^\intercal}{n}\right)\left(\mathbf{q}_{k+1}+\tilde{\gamma}\mathbf{e}_{k+1}^s\right)-\left(\mathbf{I}-\frac{\mathbf{1}\mathbf{1}^\intercal}{n}\right)\left(\mathbf{q}_{k}+\tilde{\gamma}\mathbf{e}_{k}^s\right)\\
		&=\left(\mathbf{W}_{\tilde{\gamma}}-\frac{\mathbf{1}\mathbf{1}^\intercal}{n}\right)(\mathbf{y}_{k}-\bar{\mathbf{y}}_{k})+\left(\mathbf{I}-\frac{\mathbf{1}\mathbf{1}^\intercal}{n}\right)\left(\mathbf{q}_{k+1}-\mathbf{q}_{k}+\tilde{\gamma}\mathbf{e}_{k+1}^s-\tilde{\gamma}\mathbf{e}_{k}^s\right).
	\end{align*}
	The preceding relationship further leads to
	\begin{align*}
		&\mathbb{E}\left[\left\|\mathbf{y}_{k+1}-\bar{\mathbf{y}}_{k+1}\right\|^2\big|\mathcal{F}_k^s\right]\\
		&=\mathbb{E}\left[\left\|\left(\mathbf{W}_{\tilde{\gamma}}-\frac{\mathbf{1}\mathbf{1}^\intercal}{n}\right)(\mathbf{y}_{k}-\bar{\mathbf{y}}_{k})+\left(\mathbf{I}-\frac{\mathbf{1}\mathbf{1}^\intercal}{n}\right)\left(\mathbf{q}_{k+1}-\mathbf{q}_{k}+\tilde{\gamma}\mathbf{e}_{k+1}^s-\tilde{\gamma}\mathbf{e}_{k}^s\right)\right\|^2\big|\mathcal{F}_k^s\right]\\
		&\le (1+\epsilon)\left\|\mathbf{W}_{\tilde{\gamma}}-\frac{\mathbf{1}\mathbf{1}^\intercal}{n}\right\|_2^2\left\|\mathbf{y}_{k}-\bar{\mathbf{y}}_{k}\right\|^2+2\left(1+\frac{1}{\epsilon}\right)\left\|\mathbf{I}-\frac{\mathbf{1}\mathbf{1}^\intercal}{n}\right\|_2^2\mathbb{E}\left[\left\|\mathbf{q}_{k+1}-\mathbf{q}_{k}\right\|^2\big|\mathcal{F}_k^s\right]\\
		&\quad+2\left(1+\frac{1}{\epsilon}\right)\left\|\mathbf{I}-\frac{\mathbf{1}\mathbf{1}^\intercal}{n}\right\|_2^2\mathbb{E}\left[\left\|\tilde{\gamma}\mathbf{e}_{k+1}^s-\tilde{\gamma}\mathbf{e}_{k}^s\right\|^2\big|\mathcal{F}_k^s\right] \\
		&\le (1+\epsilon)(1-\tilde{\gamma} \eta_w)^2\left\|\mathbf{y}_{k}-\bar{\mathbf{y}}_{k}\right\|^2+2\left(1+\frac{1}{\epsilon}\right)\mathbb{E}\left[\left\|\mathbf{q}_{k+1}-\mathbf{q}_{k}\right\|^2\big|\mathcal{F}_k^s\right]\\
		&\quad+2\left(1+\frac{1}{\epsilon}\right)\mathbb{E}\left[\left\|\tilde{\gamma}\mathbf{e}_{k+1}^s-\tilde{\gamma}\mathbf{e}_{k}^s\right\|^2\big|\mathcal{F}_k^s\right] ,
	\end{align*}
	where the first inequality follows from the fact $\|\mathbf{AB}\|\le \|\mathbf{A}\|_2\|\mathbf{B}\|$ for any matrices $\mathbf{A},\mathbf{B}\in \mathbb{R}^{n\times n} $ \cite[Page 2]{Song2024Optgt} and $(a+b)^2\le (1+\epsilon)a^2+\left(1+\frac{1}{\epsilon}\right)b^2$ for any scalars $a,b$ and $\epsilon>0$, the second inequality follows from the facts $\left\|\mathbf{W}_{\tilde{\gamma}}-\frac{\mathbf{1}\mathbf{1}^\intercal}{n}\right\|\le 1-\tilde{\gamma} \eta_w$ and $\left\|\mathbf{I}-\frac{\mathbf{1}\mathbf{1}^\intercal}{n}\right\|_2=1$.
	By setting $\epsilon=\frac{\tilde{\gamma} \eta_w}{1-\tilde{\gamma} \eta_w}$, we have
	\begin{align}\label{y-bound-1}
		\mathbb{E}\left[\left\|\mathbf{y}_{k+1}-\bar{\mathbf{y}}_{k+1}\right\|^2\big|\mathcal{F}_k^s\right]
		&\le (1-\tilde{\gamma} \eta_w)\left\|\mathbf{y}_{k}-\bar{\mathbf{y}}_{k}\right\|^2+\frac{2}{\tilde{\gamma} \eta_w}\mathbb{E}\left[\left\|\mathbf{q}_{k+1}-\mathbf{q}_{k}\right\|^2\big|\mathcal{F}_k^s\right]\notag\\
		&\quad+\frac{2}{\tilde{\gamma} \eta_w}\mathbb{E}\left[\left\|\tilde{\gamma}\mathbf{e}_{k+1}^s-\tilde{\gamma}\mathbf{e}_{k}^s\right\|^2\big|\mathcal{F}_k^s\right].
	\end{align}
	
	Next, we proceed to bound the last two terms on the right hand side of (\ref{y-bound-1}) separately. 
	For the term $\mathbb{E}\left[\left\|\mathbf{q}_{k+1}-\mathbf{q}_{k}\right\|^2\big|\mathcal{F}_k^s\right]$, by recursion (\ref{alg:q}), we have
	\begin{align}
		&\mathbb{E}\left[\left\|\mathbf{q}_{k+1}-\mathbf{q}_{k}\right\|^2\big|\mathcal{F}_k^s\right]\notag\\
		&=\mathbb{E}\left[\left\|-\lambda_k\mathbf{q}_{k}+\lambda_k\tilde{\mathbf{g}}_{k+1}^{'}+(1-\lambda_k)\left(\tilde{\mathbf{g}}_{k+1}^{'}-\tilde{\mathbf{g}}_k\right)
		\right\|^2\big|\mathcal{F}_k^s\right]\notag\\
		&=\mathbb{E}\left[\left\|-\lambda_k\left(\mathbf{q}_{k}-\mathbf{g}_{k}\right)+\lambda_k\left(\tilde{\mathbf{g}}_{k+1}^{'}-\mathbf{g}_{k+1}\right)+\lambda_k\left(\mathbf{g}_{k+1}-\mathbf{g}_{k}\right)+(1-\lambda_k)\left(\tilde{\mathbf{g}}_{k+1}^{'}-\tilde{\mathbf{g}}_k\right)\right\|^2\big|\mathcal{F}_k^s\right]\notag\\
		&\le 4\lambda_k^2\left\|\mathbf{e}_{k}^q\right\|^2+4\lambda_k^2\mathbb{E}\left[\left\|\tilde{\mathbf{g}}_{k+1}^{'}-\mathbf{g}_{k+1}\right\|^2\big|\mathcal{F}_k^s\right]+4\lambda_k^2\mathbb{E}\left[\left\|\mathbf{g}_{k+1}-\mathbf{g}_{k}\right\|^2\big|\mathcal{F}_k^s\right]\notag\\
		&\quad+4(1-\lambda_k)^2\mathbb{E}\left[\left\|\tilde{\mathbf{g}}_{k+1}^{'}-\tilde{\mathbf{g}}_k\right\|^2\big|\mathcal{F}_k^s\right]\notag\\
		&\le 4\lambda_k^2\left\|\mathbf{e}_{k}^q\right\|^2+4\lambda_k^2n\sigma_\xi^2+4\left(\lambda_k^2+(1-\lambda_k)^2\right)L^2\mathbb{E}\left[\left\|\mathbf{x}_{k+1}-\mathbf{x}_{k}\right\|^2\big|\mathcal{F}_k^s\right]\notag\\
		%&\le 4\lambda_k^2\mathbb{E}\left[\left\|\mathbf{e}_{k}^q\right\|^2\right]+4\lambda_k^2\sigma_\xi^2+4L^2\mathbb{E}\left[\left\|\mathbf{x}_{k+1}-\mathbf{x}_{k}\right\|^2\right]\\
		&\le 4\lambda_k^2\left\|\mathbf{e}_{k}^q\right\|^2+4\lambda_k^2n\sigma_\xi^2+20\left(1+2\lambda_k^2\right)L^2\left(\gamma_k^2\left\|\mathbf{W}-\mathbf{I}\right\|^2\left\|\mathbf{x}_k-\bar{\mathbf{x}}_k\right\|^2+\gamma_k^2\mathbb{E}\left[\left\|\mathbf{e}_k^x\right\|^2\big|\mathcal{F}_k^s\right]\right.\notag\\
		&\quad\left.+\alpha_k^2\left\|\mathbf{y}_k-\bar{\mathbf{y}}_k\right\|^2+\alpha_k^2n\left\|\bar{q}_k\right\|_2^2+\alpha_k^2\tilde{\gamma}^2\left\|\mathbf{e}_k^s\right\|^2\right),\label{g-bound-1}
	\end{align}
	where $\mathbf{e}_{k}^q=\mathbf{q}_{k}-\mathbf{g}_{k}$,  the second inequality follows from Assumption \ref{ass:stochastic gradient},  the third inequality uses the upper bound of $\mathbb{E}\left[\left\|\mathbf{x}_{k+1}-\mathbf{x}_{k}\right\|^2\big|\mathcal{F}_k^s\right]$ in (\ref{x-bound}) and the fact $\lambda_k^2+(1-\lambda_k)^2\le 1+2\lambda_k^2$. Towards the last term on the right hand side of (\ref{y-bound-1}), we have
	\begin{align}
		\mathbb{E}\left[\left\|\tilde{\gamma}\mathbf{e}_{k+1}^s-\tilde{\gamma}\mathbf{e}_{k}^s\right\|^2\big|\mathcal{F}_k^s\right]&=\tilde{\gamma}^2\mathbb{E}\left[\left\|-\beta_k\mathbf{e}_{k}^s+\beta_k\mathbf{W}\zeta_{k+1}^s\right\|^2\big|\mathcal{F}_k^s\right]\notag\\
		&\le 2\tilde{\gamma}^2\beta_k^2\left\|\mathbf{e}_{k}^s\right\|^2+2\tilde{\gamma}^2\beta_k^2\mathbb{E}\left[\left\|\mathbf{W}\zeta_{k+1}^s\right\|^2\big|\mathcal{F}_k^s\right]\notag\\
		&\le 2\tilde{\gamma}^2\beta_k^2\left\|\mathbf{e}_{k}^s\right\|^2+2\tilde{\gamma}^2\beta_k^2\sum_{i=1}^n\sum_{j\in\mathcal{N}_{i}}w_{ij
		}\mathbb{E}\left[\left\|\zeta_{j,k+1}^s\right\|^2\big|\mathcal{F}_k^s\right]\notag\\
		&\le2\tilde{\gamma}^2\beta_k^2\left\|\mathbf{e}_{k}^s\right\|^2+2\tilde{\gamma}^2\beta_k^2n\sigma_\zeta^2,\label{es-bound}
	\end{align}
	where the last inequality follows from the fact that $\sum_{i=1}^nw_{ij}= 1$ and $\mathbb{E}\left[\left\|\zeta_{j,k+1}^x\right\|^2\big|\mathcal{F}_k^s\right]=\mathbb{E}\left[\mathbb{E}\left[\left\|\zeta_{j,k+1}^x\right\|^2\big|\mathcal{F}_{k+1}^s\right]\big|\mathcal{F}_k^s\right]\le \sigma_\zeta^2$.

	Substituting (\ref{g-bound-1}) and (\ref{es-bound}) into (\ref{y-bound-1}), we arrive at
	\begin{align*}
		&\mathbb{E}\left[\left\|\mathbf{y}_{k+1}-\bar{\mathbf{y}}_{k+1}\right\|^2\big|\mathcal{F}_k^s\right]\\
		&\le \left(1-\tilde{\gamma} \eta_w+\frac{40L^2\left(1+2\lambda_k^2\right)}{\tilde{\gamma} \eta_w}\alpha_k^2\right)\left\|\mathbf{y}_{k}-\bar{\mathbf{y}}_{k}\right\|^2+\frac{40L^2\left(1+2\lambda_k^2\right)}{\tilde{\gamma} \eta_w}\left(\gamma_k^2\mathbb{E}\left[\left\|\mathbf{e}_k^x\right\|^2\big|\mathcal{F}_k^s\right]\right.\notag\\
		&\quad\left.+\gamma_k^2\left\|\mathbf{W}-\mathbf{I}\right\|^2\left\|\mathbf{x}_k-\bar{\mathbf{x}}_k\right\|^2+n\alpha_k^2\left\|\bar{q}_k\right\|_2^2\right)+\frac{8\lambda_k^2}{\tilde{\gamma} \eta_w}\left(\left\|\mathbf{e}_k^q\right\|^2+n\sigma_\xi^2\right)\notag\\
		&\quad+\frac{4\tilde{\gamma}\left(\beta_k^2+10L^2\left(1+2\lambda_k^2\right)\alpha_k^2\right)}{\eta_w}\left\|\mathbf{e}_{k}^s\right\|^2+\frac{4\tilde{\gamma}\beta_k^2}{ \eta_w}n\sigma_\zeta^2.
	\end{align*}
\end{proof}

\noindent\textbf{A.3. Proof of Lemma \ref{lem:x-1}}
\begin{proof}
	By the definition of  $\bar{\mathbf{x}}_{k+1}$,
	\begin{align*}
		\mathbf{x}_{k+1}-\bar{\mathbf{x}}_{k+1}&=\left(\mathbf{W}_{\gamma_k}-\frac{\mathbf{1}\mathbf{1}^\intercal}{n}\right)\mathbf{x}_{k}-\alpha_k\left(\mathbf{y}_{k}-\bar{\mathbf{y}}_{k}\right)+\gamma_k\left(\mathbf{I}-\frac{\mathbf{1}\mathbf{1}^\intercal}{n}\right)\mathbf{e}_{k}^x\\
		&=\left(\mathbf{W}_{\gamma_k}-\frac{\mathbf{1}\mathbf{1}^\intercal}{n}\right)(\mathbf{x}_{k}-\bar{\mathbf{x}}_{k})-\alpha_k\left(\mathbf{y}_{k}-\bar{\mathbf{y}}_{k}\right)+\gamma_k\left(\mathbf{I}-\frac{\mathbf{1}\mathbf{1}^\intercal}{n}\right)\mathbf{e}_{k}^x.
	\end{align*}
	Taking square norms and conditional expectation on both sides of the above equality,
	\begin{align}
		&\mathbb{E}\left[\left\|\mathbf{x}_{k+1}-\bar{\mathbf{x}}_{k+1}\right\|^2\big|\mathcal{F}_k^s\right]\notag\\
		&=\mathbb{E}\left[\left\|\left(\mathbf{W}_{\gamma_k}-\frac{\mathbf{1}\mathbf{1}^\intercal}{n}\right)(\mathbf{x}_{k}-\bar{\mathbf{x}}_{k})-\alpha_k\left(\mathbf{y}_{k}-\bar{\mathbf{y}}_{k}\right)+\gamma_k\left(\mathbf{I}-\frac{\mathbf{1}\mathbf{1}^\intercal}{n}\right)\mathbf{e}_{k}^x\right\|^2\big|\mathcal{F}_k^s\right]\notag\\
		&\le (1+\epsilon)\left\|\mathbf{W}_{\gamma_k}-\frac{\mathbf{1}\mathbf{1}^\intercal}{n}\right\|_2^2\left\|\mathbf{x}_{k}-\bar{\mathbf{x}}_{k}\right\|^2+\left(1+\frac{1}{\epsilon}\right)\mathbb{E}\left[\left\|-\alpha_k\left(\mathbf{y}_{k}-\bar{\mathbf{y}}_{k}\right)+\gamma_k\left(\mathbf{I}-\frac{\mathbf{1}\mathbf{1}^\intercal}{n}\right)\mathbf{e}_{k}^x\right\|^2\big|\mathcal{F}_k^s\right]\notag\\
		&\le (1+\epsilon)(1-\gamma_k \eta_w)^2\left\|\mathbf{x}_{k}-\bar{\mathbf{x}}_{k}\right\|^2+2\alpha_k^2\left(1+\frac{1}{\epsilon}\right)\left\|\mathbf{y}_{k}-\bar{\mathbf{y}}_{k}\right\|^2+2\gamma_k^2\left(1+\frac{1}{\epsilon}\right)\mathbb{E}\left[\left\|\mathbf{e}_{k}^x\right\|^2\big|\mathcal{F}_k^s\right]\notag\\
		&\le (1-\gamma_k \eta_w)\left\|\mathbf{x}_{k}-\bar{\mathbf{x}}_{k}\right\|^2+\frac{2\alpha_k^2}{\gamma_k \eta_w}\left\|\mathbf{y}_{k}-\bar{\mathbf{y}}_{k}\right\|^2+\frac{2\gamma_k}{ \eta_w}\mathbb{E}\left[\left\|\mathbf{e}_{k}^x\right\|^2\big|\mathcal{F}_k^s\right],
	\end{align}
	where $\eta_w$ is defined in (\ref{eta}),  the first inequality follows from the fact $\|\mathbf{AB}\|\le \|\mathbf{A}\|_2\|\mathbf{B}\|$ for any matrices $\mathbf{A},\mathbf{B}\in \mathbb{R}^{n\times n} $ and $(a+b)^2\le (1+\epsilon)a^2+\left(1+\frac{1}{\epsilon}\right)b^2$ for any scalars $a,b$ and $\epsilon>0$, the second inequality follows from the facts $\left\|\mathbf{W}_{\gamma_k}-\frac{\mathbf{1}\mathbf{1}^\intercal}{n}\right\|\le 1-\gamma_k \eta_w$ and  $\left\|\mathbf{I}-\frac{\mathbf{1}\mathbf{1}^\intercal}{n}\right\|_2=1$, the third inequality follows from  the setting $\epsilon=\frac{\gamma_k \eta_w}{1-\gamma_k \eta_w}$.
\end{proof}

\noindent\textbf{A.4. Proof of Lemma \ref{lem:xstar-1}}

\begin{proof}
	A similar result appears in \cite[Lemma 2]{xin21Hybrid}. For completeness, the full proof is included here. By the Lipschitz smoothness of $f(x)$,
	\begin{align*}
		f(\bar{x}_{k+1})&\le f(\bar{x}_k)-\langle  \nabla f(\bar{x}_k), \bar{x}_{k+1}-\bar{x}_k\rangle+\frac{L}{2}\|\bar{x}_{k+1}-\bar{x}_k\|_2^2\\
		&= f(\bar{x}_k)-\alpha_k\langle  \nabla f(\bar{x}_k), (\bar{x}_{k+1}-\bar{x}_k)/\alpha_k\rangle+\frac{L}{2}\|\bar{x}_{k+1}-\bar{x}_k\|_2^2\\
		&=f(\bar{x}_k)-\frac{\alpha_k}{2}\|\nabla f(\bar{x}_k)\|_2^2-\left(\frac{2}{\alpha_k}-\frac{L}{2}\right)\|\bar{x}_{k+1}-\bar{x}_k\|_2^2+\frac{2}{\alpha_k}\left\|\bar{x}_{k+1}-\bar{x}_k+\alpha_k\nabla f(\bar{x}_k)\right\|_2^2.
	\end{align*}
	Note that $\bar{x}_{k+1}-\bar{x}_k=-\alpha_k\bar{y}_{k}+\gamma_k\bar{e}^x_{k}=-\alpha_k(\bar{q}_{k}+\tilde{\gamma}\bar{e}^s_k)+\gamma_k\bar{e}^x_{k}$ and $\alpha_k\le \frac{1}{L}$, 
	%\begin{align*}
	%	\|\bar{x}_{k+1}-\bar{x}_k\|_2^2=\|-\alpha_k\bar{y}_{k}+\gamma_k\bar{e}^x_{k}\|^2&=\|-\alpha_k(\bar{q}_{k}+\tilde{\gamma}\bar{e}^s_k)+\gamma_k\bar{e}^x_{k}\|_2^2\\
	%	&\le 3\alpha_k^2\|\bar{q}_{k}\|_2^2+3\alpha_k^2\|\bar{e}^s_k\|^2+3\gamma_k^2\|\bar{e}^x_{k}\|^2\\
	%	&\le 3\alpha_k^2\|\bar{q}_{k}\|_2^2+\frac{3\alpha_k^2}{n}\|\mathbf{e}^s_k\|^2+\frac{3\gamma_k^2}{n}\|\mathbf{e}^x_{k}\|^2,
	%\end{align*}
	we have 
	\begin{align}\label{f-bound}
		f(\bar{x}_{k+1})
		&\le f(\bar{x}_k)-\frac{\alpha_k}{2}\|\nabla f(\bar{x}_k)\|_2^2-\frac{1}{\alpha_k}\|-\alpha_k(\bar{q}_{k}+\tilde{\gamma}\bar{e}^s_k)+\gamma_k\bar{e}^x_{k}\|_2^2\notag\\
		&\quad+\frac{2}{\alpha_k}\left\|-\alpha_k(\bar{q}_{k}+\tilde{\gamma}\bar{e}^s_k)+\gamma_k\bar{e}^x_{k}+\alpha_k\nabla f(\bar{x}_k)\right\|_2^2\notag\\
		&\le  f(\bar{x}_k)-\frac{\alpha_k}{2}\|\nabla f(\bar{x}_k)\|_2^2-\frac{1}{\alpha_k}\|-\alpha_k(\bar{q}_{k}+\tilde{\gamma}\bar{e}^s_k)+\gamma_k\bar{e}^x_{k}\|_2^2\notag\\
		&\quad+4\alpha_k\left(\|\bar{q}_{k}-\nabla f(\bar{x}_k)\|_2^2+\left\|\tilde{\gamma}\bar{e}^s_k+\frac{\gamma_k}{\alpha_k}\bar{e}^x_{k}\right\|_2^2\right)\notag\\
		&=f(\bar{x}_k)-\frac{\alpha_k}{2}\|\nabla f(\bar{x}_k)\|_2^2-\alpha_k\|\bar{q}_{k}\|_2^2-\alpha_k\left\|\tilde{\gamma}\bar{e}^s_k+\frac{\gamma_k}{\alpha_k}\bar{e}^x_{k}\right\|_2^2+2\alpha_k\left\langle \bar{q}_{k},\tilde{\gamma}\bar{e}^s_k+\frac{\gamma_k}{\alpha_k}\bar{e}^x_{k}\right\rangle\notag\\
		&\quad+4\alpha_k\left(\|\bar{q}_{k}-\nabla f(\bar{x}_k)\|_2^2+\left\|\tilde{\gamma}\bar{e}^s_k+\frac{\gamma_k}{\alpha_k}\bar{e}^x_{k}\right\|_2^2\right)\notag\\
		&\le f(\bar{x}_k)-\frac{\alpha_k}{2}\|\nabla f(\bar{x}_k)\|_2^2-\alpha_k\|\bar{q}_{k}\|_2^2+4\alpha_k\|\bar{q}_{k}-\nabla f(\bar{x}_k)\|_2^2+5\alpha_k\left\|\tilde{\gamma}\bar{e}^s_k+\frac{\gamma_k}{\alpha_k}\bar{e}^x_{k}\right\|_2^2
	\end{align}
	where the second and third inequalities use Young’s and Cauchy-Schwarz inequalities. Toward the last two terms on the right hand side of the above inequality, we have 
	\begin{align*}
		&5\alpha_k\left\|\tilde{\gamma}\bar{e}^s_k+\frac{\gamma_k}{\alpha_k}\bar{e}^x_{k}\right\|_2^2\le \frac{10\alpha_k\tilde{\gamma}^2}{n}\|\mathbf{e}^s_{k}\|^2+\frac{10\gamma_k^2}{\alpha_kn}\|\mathbf{e}^x_{k}\|^2,
	\end{align*}
	and
	\begin{align*}
		4\alpha_k\|\bar{q}_{k}-\nabla f(\bar{x}_k)\|_2^2=&4\alpha_k\left\|\bar{q}_{k}-\frac{1}{n}\sum_{j=1}^n\nabla f_j(x_{j,k})+\frac{1}{n}\sum_{j=1}^n\nabla f_j(x_{j,k})-\nabla f(\bar{x}_k)\right\|_2^2\\
		\le&8\alpha_k\left\|\bar{q}_{k}-\frac{1}{n}\sum_{j=1}^n\nabla f_j(x_{j,k})\right\|_2^2+8\alpha_k\left\|\frac{1}{n}\sum_{j=1}^n\nabla f_j(x_{j,k})-\nabla f(\bar{x}_k)\right\|_2^2\\
		\le&\frac{8\alpha_k}{n}\left\|\mathbf{e}^q_k\right\|^2+\frac{8\alpha_kL^2}{n}\left\|\mathbf{x}_k-\bar{\mathbf{x}}_k\right\|^2.
	\end{align*}
	Therefore,
	\begin{align*}
		f(\bar{x}_{k+1})
		&\le f(\bar{x}_k)-\frac{\alpha_k}{2}\|\nabla f(\bar{x}_k)\|_2^2-\alpha_k\|\bar{q}_{k}\|_2^2+\frac{8\alpha_k}{n}\left\|\mathbf{e}^q_k\right\|^2+\frac{8\alpha_kL^2}{n}\left\|\mathbf{x}_k-\bar{\mathbf{x}}_k\right\|^2\\
		&\quad+\frac{10\alpha_k\tilde{\gamma}^2}{n}\|\mathbf{e}^s_{k}\|^2+\frac{10\gamma_k^2}{\alpha_kn}\|\mathbf{e}^x_{k}\|^2.
	\end{align*}
	Taking the conditional expectation on both sides of the above inequality,  then by the the strong convexity of $f(x)$, we arrive at
	\begin{align*}
		\mathbb{E}\left[f(\bar{x}_{k+1})-f^*\big|\mathcal{F}_k^s\right]	
		&\le (1-\mu\alpha_k)(f(\bar{x}_k)-f^*)-\alpha_k\|\bar{q}_{k}\|_2^2+\frac{8\alpha_k}{n}\left\|\mathbf{e}^q_k\right\|^2+\frac{8\alpha_kL^2}{n}\left\|\mathbf{x}_k-\bar{\mathbf{x}}_k\right\|^2\\
		&\quad+\frac{10\alpha_k\tilde{\gamma}^2}{n}\|\mathbf{e}^s_{k}\|^2+\frac{10\gamma_k^2}{\alpha_kn}\mathbb{E}\left[\|\mathbf{e}^x_{k}\|^2\big|\mathcal{F}_k^s\right],
	\end{align*}
	where $f^*=\min_{x\in\mathbb{R}^d} f(x)$.
\end{proof}

\section*{Appendix B: Technical Results}
The following two lemmas are technical results for establishing the almost sure convergence rate of VRA-D(S)GT. The first lemma shows that $\mathcal{L}^S_k$ is finite almost surely, while the second lemma, adapted from \cite[Theorem F.1]{Peggy2025Newton}, establishes almost sure convergence of the martingale sequence.
\begin{lem}\label{lem:e-o}
Under the conditions of Theorem \ref{thm:as-DSGT}, $\mathcal{L}^S_k$ converges to a finite random variable almost surely.
\end{lem}
\begin{proof}
	 By (\ref{L-bound-4}) and Proposition \ref{prop:vra-1},
	\begin{align}\label{L-bound-5}
		&\mathbb{E}\left[\mathcal{L}^S_{k+1}\big|\mathcal{F}_k^s\right]\notag\\
		&\le\mathcal{L}^S_{k}+\left[\left(\frac{80L^2}{\tilde{\gamma} \eta_w}+30L^2\right)\gamma_k^2+\frac{2\gamma_k}{ \eta_w}+\frac{10\gamma_k^2}{\alpha_kn}\right]\mathbb{E}\left[\left\|\mathbf{e}_k^x\right\|^2\big|\mathcal{F}_k^s\right] +\frac{8\lambda_k^2}{\tilde{\gamma} \eta_w}\left(n\sigma_\xi^2\right)+\frac{4\tilde{\gamma}n\sigma_\zeta^2\beta_k^2}{ \eta_w}\notag\\
		&\quad+3n\sigma_\xi^2\lambda_k^2+\left[150L^2\alpha_k^2\tilde{\gamma}^2+\frac{4\tilde{\gamma}\left(\beta_k^2+20L^2\alpha_k^2\right)}{\eta_w}+\frac{10\alpha_k\tilde{\gamma}^2}{n}\right]\tilde{c}_s\beta_k\ln (k+1)\notag\\
%		&\le\mathcal{L}^S_{k}+\left[\left(\frac{80L^2}{\tilde{\gamma} \eta_w}+30L^2\right)c_1^2\alpha_1+\frac{2c_1}{ \eta_w}+\frac{10c_1^2}{n}\right]\alpha_{k}\mathbb{E}\left[\left\|\mathbf{e}_k^x\right\|^2\big|\mathcal{F}_k^s\right] +\frac{8c_2^2\alpha_k^2}{\tilde{\gamma} \eta_w}\left(n\sigma_\xi^2\right)+\frac{4\tilde{\gamma}n\sigma_\zeta^2c_3^2\alpha_k^2}{ \eta_w}\notag\\
%		&\quad+3n\sigma_\xi^2c_2^2\alpha_k^2+\left[150L^2\alpha_1\tilde{\gamma}^2+\frac{4\tilde{\gamma}\left(c_3^2\alpha_1+20L^2\alpha_1\right)}{\eta_w}+\frac{10\tilde{\gamma}^2}{n}\right]\tilde{c}_sc_3\alpha_k^2\ln (k+1)\\
		&\le\mathcal{L}^S_{k}+C_1\alpha_{k}\mathbb{E}\left[\left\|\mathbf{e}_k^x\right\|^2\big|\mathcal{F}_k^s\right] +C_2\alpha_k^2\ln (k+1),
	\end{align}
	where 
	\begin{align*}
	&C_1=\left[\left(\frac{80L^2}{\tilde{\gamma} \eta_w}+30L^2\right)c_1^2\alpha_1+\frac{2c_1}{ \eta_w}+\frac{10c_1^2}{n}\right],\\
	&C_2=\frac{8c_2^2}{\tilde{\gamma} \eta_w}\left(n\sigma_\xi^2\right)+\frac{4\tilde{\gamma}n\sigma_\zeta^2c_3^2}{ \eta_w}+3n\sigma_\xi^2c_2^2+\left[150L^2\alpha_1\tilde{\gamma}^2+\frac{4\tilde{\gamma}\left(c_3^2\alpha_1+20L^2\alpha_1\right)}{\eta_w}+\frac{10\tilde{\gamma}^2}{n}\right]\tilde{c}_sc_3,
	\end{align*}
	the second inequality uses the fact that $\gamma_k=c_1\alpha_k,\lambda_k=c_2\alpha_k$ and $\beta_k\le c_3 \alpha_k$. Note that the last term on the right-hand side of (\ref{L-bound-5}) is summable, that is,
	$$\sum_{k=1}^\infty C_2\alpha_k^2\ln (k+1)<\infty.$$
Then  Robbins–Siegmund supermartingale convergence theorem \cite{Robbins1971A} implies that $\mathcal{L}^S_k$ converges to a finite random variable almost surely if the second term on the right-hand side of (\ref{L-bound-5}) satisfies
\begin{equation}\label{summable}
	\sum_{k=1}^\infty C_1\alpha_{k}\mathbb{E}\left[\left\|\mathbf{e}_k^x\right\|^2\big|\mathcal{F}_k^s\right]<\infty~~~\text{a.s.}.
\end{equation}
 
We now show (\ref{summable}). According to the definition of $\mathbf{e}^x_{k}$ in (\ref{e}) and recursion of $z_{i,k+1}^x$ in (\ref{alg:z-x}),
	\begin{equation*}
	\mathbf{e}^x_{k}	=(1-\beta_{k-1})\mathbf{e}^x_{k-1}+\beta_{k-1} \tilde{\mathbf{W}}\zeta_{k}^x.
	\end{equation*}
	Thus,
	\begin{align}\label{e-bound}
	\mathbb{E}\left[\left\|\mathbf{e}_k^x\right\|^2\big|\mathcal{F}_k^s\right]
	&=(1-\beta_{k-1})^2\left\|\mathbf{e}_{k-1}^x\right\|^2+\beta_{k-1}^2\mathbb{E}\left[\left\|\tilde{\mathbf{W}}\zeta_{k}^x\right\|^2\big|\mathcal{F}_k^s\right]\notag\\
	&\le (1-\beta_{k-1})^2\left\|\mathbf{e}_{k-1}^x\right\|^2+n\sigma_\zeta^2\beta_{k-1}^2\notag\\
	&\le (1-\beta_{k-1})^2\tilde{c}_x\beta_{k-1}\ln (k)+n\sigma_\zeta^2\beta_{k-1}^2~~~\text{a.s.}
	\end{align}
	where $\tilde{c}_x$ is some positive finite random variable, the equality follows from the condition  $\mathbb{E}\left[\zeta_{i,k}^x\big|\mathcal{F}_k^s\right]=\mathbf{0}$,  the first inequality follows from the condition $ \mathbb{E}\left[\|\zeta^x_{i,k}\|^2\big|\mathcal{F}_k^s\right]\le\sigma_\zeta^2$, the second inequality follows from Proposition \ref{prop:vra-1}.
Therefore, 
	\begin{align*}
		\sum_{k=1}^\infty C_1\alpha_{k}\mathbb{E}\left[\left\|\mathbf{e}_k^x\right\|^2\big|\mathcal{F}_k^s\right]&\le C_1\alpha_1\left\|\mathbf{e}_1^x\right\|^2+\sum_{k=2}^\infty C_1\alpha_k\left[(1-\beta_{k-1})^2\tilde{c}_x\beta_{k-1}\ln (k)+n\sigma_\zeta^2\beta_{k-1}^2\right]\\
		&\le C_1\alpha_1\left\|\mathbf{e}_1^x\right\|^2+ \max_k\left\{(1-\beta_{k})^2\tilde{c}_x+n\sigma_\zeta^2\beta_{k}\right\}\sum_{k=1}^\infty C_1c_3\alpha_k^2\ln (k+1)\\
		&<\infty~~~\text{a.s.}
	\end{align*}
 where the second inequality follows from the fact $\beta_k\le c_3 \alpha_k$ and that $\alpha_{k}$ is nonincreasing. The proof is complete.
\end{proof}

\begin{lem}\label{lem:as}
Let $(\Omega,\mathcal{F},\P)$ be a probability space and  
\begin{equation}\label{M-def}
	M_{k+1}=\sum_{t=1}^kB_{k,t}\beta_{t}R_t\tilde{\zeta}_{t+1},
\end{equation}
where
\begin{itemize}
	\item[(C1)] $\{\tilde{\zeta}_k\}$ is a martingale difference sequence in $\mathbb{R}^{n\times d}$ adapted to a filtration $\{\mathcal{F}_k\}$ such that
	\begin{align*}
	&\mathbb{E}\left[\|\tilde{\zeta}_{k+1}\|^2 \big|\mathcal{F}_k\right]\le C+R_{2,k}~a.s.,\\
	&\sum_{t= 1}^\infty \beta_t\mathbb{E}\left[\|\tilde{\zeta}_{t+1}\|^21_{\left\{\|\tilde{\zeta}_{k+1}\|^2\ge \frac{1}{\beta_k\ln (k+1)}\right\}}\big|\mathcal{F}_k\right]< \infty;
	\end{align*}
	where $C\ge 0$ and $\{R_{2,k}\}$ converges almost surely to 0;
	\item[(C2)] $\beta_{k}\asymp\frac{1}{k^{b_1}}$ with  and $b_1\in(0.5,1)$;
	\item[(C3)] $\{R_{k}\}$ is a sequence  of matrices such that, for a deterministic  sequence $\{v_{k}\}$,
	\begin{equation*}
	\|R_k\|=o(v_k)~\text{or}~\|R_k\|=\mathcal{O}(v_k)
	\end{equation*}
	where $v_k=\frac{(\ln (k))^{b_3}}{k^{b_2}}$ with $b_3,b_2\ge 0$;
	\item[(C4)]  $$B_{k,t}=\Pi_{l=t+1}^k(\mathbf{I}-\beta_{l}\Gamma),~B_{k,k+1}=\mathbf{I},$$
	 or $$B_{k,t}=\tilde{\Gamma}^{k-t}$$ for all $k\ge 1$ and $1\le t\le k$, where $\Gamma\in \mathbb{R}^{n\times n}$ is a positive definite matrix, $\tilde{\Gamma}\in \mathbb{R}^{n\times n}$ satisfies $\|\tilde{\Gamma}\|_2<1$.
\end{itemize}
Then,
\begin{equation*}\label{M-rate}
\|M_k\|=\mathcal{O}\left(\beta_{k}v_k^2\ln (k+1)\right).
\end{equation*}
\end{lem}
%\begin{proof}
%When $B_{k,t}$ is defined as in (\ref{B-1}), the convergence result (\ref{M-rate}) follows directly from \cite[Theorem F.1]{Peggy2025Newton}. Note that for all sufficiently large $k$,
%\begin{equation*}
%\left\|\tilde{\Gamma}\right\|_2\le \left\|\mathbf{I}-\beta_k\Gamma\right\|_2,~\left\|\tilde{\Gamma}\right\|_2^{k-t+1}\le c \left\|\Pi_{l=t+1}^k(\mathbf{I}-\beta_{l}\Gamma)\right\|_2,
%\end{equation*}
%for some positive constant $c$.
%By extending the proof of \cite[Theorem F.1]{Peggy2025Newton}, we can further show that (\ref{M-rate}) remains valid when $B_{k,t}$ is given by (\ref{B-2}).
%\end{proof}

\end{document}